\documentclass[10pt]{article}
\usepackage[multiple]{footmisc}
\usepackage{pgf,tikz,pgfplots}
\pgfplotsset{compat=1.15}
\usepackage{mathrsfs}
\usepackage[utf8]{inputenc}
\usepackage[T2A]{fontenc}
\usepackage[russian,english]{babel}
\usepackage{hyperref}

\usepackage{amsthm, amsmath}
\usepackage[shortlabels]{enumitem}

\usetikzlibrary{arrows}

\usepackage{fullpage}
\usepackage{graphicx}
\usepackage{subfigure} 

\usepackage[h]{esvect}

\newtheorem{theorem}{Theorem}[section]

\newtheorem{definition}[theorem]{Definition}
\newtheorem{proposition}[theorem]{Proposition}
\newtheorem{corollary}[theorem]{Corollary}
\newtheorem{remark}[theorem]{Remark}
\newtheorem{lemma}[theorem]{Lemma}
\newtheorem{example}[theorem]{Example}

\usepackage{amsmath}
\usepackage{amssymb}
        \def\X{\mathrm{X}}
        \def\E{\mathrm{E}}
        \def\S{\mathrm{S}}
        \def\G{\mathrm{G}}
        \def\F{\mathrm{F}}
\usepackage{url}
\usepackage{graphicx}

\newcommand\N{\mathbb{N}}
\newcommand\R{\mathbb{R}}

\newcommand{\St}{{\mathcal St}}

        \newcommand{\HH}{{\mathcal H}}
        \renewcommand{\H}{\HH^1}

   \newcommand{\defeq}{:=}
        \def\dist{\mathrm{dist}\,}
        \def\conv{\mathrm{conv}\,}
        \def\oord{\mathrm{ord}\,}
        
        \def\oordball{\mathrm{ordball}\,}
        \def\Int{\mathrm{Int}\,}
        
        \def\diam{\mathrm{diam}\,}

\usepackage[textwidth=2cm]{todonotes}

\title{On regularity of maximal distance minimizers in $\mathbb{R}^n$}

\author{Alexey Gordeev\footnote{The Euler International Mathematical Institute, Fontanka 27, St. Petersburg, 191023, Russia.}\, and Yana Teplitskaya\footnote{Chebyshev Laboratory, St. Petersburg State University, 14th Line V.O., 29, Saint Petersburg 199178 Russia.}\footnote{Université Paris-Saclay, rue Michel Magat Building 307, Orsay Cedex, 91405, France.}}

\begin{document}

\maketitle

\renewcommand{\abstractname}{Abstract}
\begin{abstract}
We study the properties of sets $\Sigma$ which are the solutions of the maximal distance minimizer problem, i.e. of sets
having the minimal length (one-dimensional Hausdorff measure) over the class	of closed connected sets $\Sigma \subset \mathbb{R}^n$ satisfying the inequality
\[
\max_{y \in M} \dist(y,\Sigma) \leq r
\]
for a given compact
set $M \subset \mathbb{R}^n$ and some given $r > 0$. 
Such sets can be considered as the shortest networks of radiating Wi-Fi cables arriving to each customer (for the set $M$ of customers) at a distance at most $r$. 

In this paper we prove that any maximal distance minimizer $\Sigma \subset \mathbb{R}^n$ has at most $3$ tangent rays at each point and the angle between any two tangent rays at the same point is at least $2\pi/3$. 
Moreover, in the plane (for $n=2$) we show that the number of points with three tangent rays is finite and every maximal distance minimizer is a finite union of simple curves with one-sided tangents continuous from the corresponding side.

 




All the results are proved for the more general class of local minimizers, i.e. sets which are optimal under a perturbation of a neighbourhood of their arbitrary point.
\end{abstract}

\section{Introduction}

For a given compact set $M \subset \mathbb{R}^n$ consider the functional
\[
	\F_{M}(\Sigma)\defeq \sup _{y\in M}\dist (y, \Sigma),
\]
where $\Sigma$ is a closed subset of $\R^n$ and $\dist(y, \Sigma)$ stands for the Euclidean distance between $y$ and $\Sigma$ (naturally, $\F_{M} (\emptyset) \defeq +\infty$). The quantity $\F_M (\Sigma)$ will be called the \emph{energy} of $\Sigma$.
Consider the class of closed connected sets $\Sigma \subset \R^n$ satisfying $\F_M(\Sigma) \leq r$ for some $r > 0$. We are interested in the properties of
sets of the minimal length (one-dimensional Hausdorff measure) $\H(\Sigma)$ over the mentioned class. Such sets will be further called \emph{minimizers}.


Let $B_\rho (x)$ be the open ball of radius $\rho$ centered at a point $x$, and let $B_\rho(M)$ be the open $\rho$-neighborhood of $M$, i.e.\
	\[
        		B_\rho(M) \defeq \bigcup_{x\in M} B_\rho(x).
	\]
	
It is easy to see that any minimizer $\Sigma$ is bounded (and, thus, compact), since $\Sigma \subset \overline{B_r(\conv M)}$, where $\conv M$ stands for the convex hull of the set $M$. 
 



In~\cite{mir} (for the plane) and in~\cite{PaoSte04max} (for $\mathbb{R}^n$) the following properties of minimizers have been proven:
    \begin{enumerate}[(a)]
        \item For all $r > 0$ the set  of minimizers is non-empty.
        \label{u}
        \item A minimizer cannot contain loops (homeomorphic images of circles).
        \label{aaaa}
        \item Any minimizer $\Sigma$ is \emph{Ahlfors regular}, i.e. there exist constants $c>0$ and $C>0$ such that for a sufficiently small $\rho>0$ for every point $x \in \Sigma$ the inequality
        \[
        c\rho \leq \H(\Sigma \cap B_\rho(x))\leq C\rho
        \]
        \label{j}
        holds.     
        \label{prop:ahlfors}
    \end{enumerate}
Also, since any minimizer $\Sigma$ is a connected closed set of a finite length, the following statement holds: 
\begin{itemize}
\item[(d)]
Any minimizer $\Sigma$ has a tangent line at almost every point of $\Sigma$. 
\label{e}

\end{itemize}
In this work (see Lemma~\ref{lm:tech} and Corollary~\ref{cor:struct}) we refine the property~\ref{prop:ahlfors}: we show that for any minimizer $\Sigma$ and any point $x \in \Sigma$ we have
\[
   \H(\Sigma \cap B_\rho(x))=\oord_x\Sigma \cdot \rho+o(\rho),
\]
where $\oord_x\Sigma \in \{1, 2, 3\}$ is defined as follows.


\begin{definition}
The order $\oord_x \Gamma$ of a point $x \in \R^n$ in a set $\Gamma \subset \mathbb{R}^n$ is the supremum of $m \in \N$ for which there exists $\rho = \rho(x) > 0$ such that for any open set $U \ni x$ with $ \diam U < \rho $ the inequality $ \sharp (\Gamma \cap \partial U) \geq m $ holds (where $\sharp A$ is the number of elements in a set $A$).
In other words, 
\[
\oord_x \Gamma := \sup\{m \in \mathbb{N}\ |\ \exists \rho>0: \sharp(\Gamma \cap \partial U)\geq m\ \forall U \text{ open }, \diam U <\rho, x \in U \}.
\]
\end{definition}
\begin{remark}
Let the order of a point $x \in \R^n$ in a set $\Gamma \subset \mathbb{R}^n$ be finite: $m \defeq \oord_x \Gamma < \infty$. Then there exists $\rho = \rho(x) > 0$ such that for any open set $U \ni x$ with $ \diam U < \rho $ the inequality $ \sharp (\Gamma \cap \partial U) \geq m $ holds, and for any positive $\varepsilon < \rho$ there exists an open set $ V \ni x $ with $ \diam V < \varepsilon $ and $ \sharp (\Gamma \cap \partial V) = m $.
\end{remark}


The main result of this article is the following theorem (see Definitions~\ref{def:tgray} and~\ref{def:one-sided} for the notions of tangent rays and one-sided tangents).
It is a special case of a broader Theorem~\ref{main-local}, stated below.

\begin{theorem}\label{mainT1}
Let $\Sigma$ be a minimizer for a compact set $M \subset \mathbb{R}^n$ and $r > 0$.
Then there are at most three tangent rays at any point of $\Sigma$,
and the pairwise angles between the tangent rays are at least $2\pi / 3$. Furthermore, tangent rays coincide with one-sided tangents, particularly the angles between one-sided tangents cannot be equal to $0$, i.e. there is one to one correspondence between tangent rays at an arbitrary point $x\in \Sigma$ and connected components of $\Sigma \setminus \{x\}.$
Moreover, if $n = 2$, then $\Sigma$ is a finite union of simple curves with one-sided tangents continuous from the corresponding side.
\end{theorem}

It is known that (see, for example, Theorem 3.1 in~\cite{thomassen2008graph}) the following lemma holds:

\begin{lemma}\label{ord_arcs}
Let $\Sigma$ be a minimizer for a compact set $M \subset \mathbb{R}^n$ and $r > 0$, and let $ x \in \Sigma $ be a point of order $ n $: $\oord_x \Sigma = n$.
Then there exists $\rho = \rho(x) > 0$ such that for any positive $\varepsilon < \rho$ the set $\Sigma \cap B_\varepsilon (x)$ contains $n$ arcs with one end at $\partial B_\varepsilon(x)$ and another end at $x$, such that any pair of these arcs intersects only at point $x$.
\end{lemma}

Note that Theorem~\ref{mainT1} implies that in the planar case the set $\Sigma \cap B_\varepsilon (x)$ coincides with these arcs.

\begin{definition}
Let $M \subset \mathbb{R}^n$ be a compact set and let $r > 0$.
A closed connected set $\Sigma \subset \mathbb{R}^n$ with $\H(\Sigma)<\infty$ is called a local minimizer if $\F_M(\Sigma) \leq r$
and there exists $\varepsilon > 0$ such that for any connected set $\Sigma'$ satisfying $\F_M(\Sigma') \leq r$
and $\diam (\Sigma \triangle \Sigma') \leq \varepsilon$ the inequality $\H(\Sigma) \leq \H(\Sigma')$ holds, where $\triangle$ is the symmetric difference.
\end{definition}

Any minimizer is also a local minimizer, so Theorem~\ref{mainT1} is a direct corollary of the following theorem.


\begin{theorem}\label{main-local}
Let $\Sigma$ be a local minimizer for a compact set $M \subset \mathbb{R}^n$ and $r > 0$.
Then $\Sigma$ has the same properties as in Theorem~\ref{mainT1}.
\end{theorem}
\begin{proof}
See Corollary~\ref{cor:struct}, Theorem~\ref{thm:plane} and Lemma~\ref{zapred2}.
\end{proof}



We should also mention the problem of minimizing the average distance functional (for more details about this problem, see the survey~\cite {L} and~\cite{lemenant2011regularity}). Maximum and average distance minimizers have a number of similar geometric and topological properties (see~\cite{PaoSte04max}). 
For example, it was proved (see~\cite{BS1}) that any minimizer of the average distance in $\mathbb{R}^2$ is a finite plane graph.
Theorem~\ref{mainT1} implies that any maximal distance minimizer in $\mathbb{R}^2$ is also a finite plane graph.
Note that this is not true for the Steiner problem (which is defined in Subsection~\ref{stp}): an example of a Steiner tree which cannot be represented as a finite plane graph can be found in~\cite{paolini2015example} (see also~\cite{paoliniSteinerTreeConnecting2023} and~\cite{cherkashinSelfSimilarInfiniteBinary2023}).

The article is organized as follows. Subsection~\ref{stp} contains the statement of the Steiner problem and some technical statements about Steiner trees (an advanced reader can skip this subsection).
Section~\ref{not} introduces notation, basic definitions and some elementary statements.
Section~\ref{res} contains the proof of the part of Theorem~\ref{main-local} about tangent rays.
Subsection~\ref{subsec:tech} contains several technical results which are used in the next section.
Finally, Section~\ref{planarcase} deals with the planar part of Theorem~\ref{main-local}.

This work is a revised version of the work~\cite{teplitskaya2019regularity}. Although the finiteness of the number of branching points (i.e. points with three tangent rays) in space remains an open problem,
the proof in the planar case follows from statements, most of which are shown to also be true in space (see Subsection~\ref{subsec:tech}).
This provides a starting point for either proving or refuting the spatial claim.



\subsection{Steiner problem} \label{stp}

Let $C \subset \mathbb{R}^n$ be a compact set (in particular, often $C=\{A_1,\ldots, A_m\}$ is supposed to be a finite set of points).
Then the shortest (having minimal one-dimensional Hausdorff measure) closed set $S$ such that $S\cup C$ is a connected set (if $C$ is a finite set of points, then one can write $C \subset S$) is called a \textit{Steiner tree}. It is known that a Steiner tree $S$ exists (but may be not unique), coincides with a (finite if $C$ is finite) union of straight segments and does not have cycles. 

Let $C=\{A_1,\ldots, A_m\}$ be a finite set. Its Steiner tree $S$ is a finite connected acyclic graph, i.e. an embedded tree. When speaking of this tree, without loss of generality, we always consider the graph with the minimal possible vertex set. It is known that vertices of this graph cannot have degree greater than $ 3 $, and that only the vertices $ A_i $ can have degree $ 1 $ and $ 2 $. Therefore, all other vertices have degree $ 3 $ and are called \textit{Steiner points} or \textit{branching points}. It follows that the number of Steiner points does not exceed $ m-2 $.

The pairwise angles between edges incident to an arbitrary vertex are at least $ 2 \pi / 3 $. Thus, a Steiner point has exactly three edges with pairwise angles $ 2 \pi / 3 $. If all vertices $ A_i $ have degree $ 1 $, then the number of Steiner points is equal to $ m-2 $ (the converse is also true), and $ S $ is called a \textit{full} Steiner tree. A set, each connected component of which is a (full) Steiner tree, is called a \textit{(full) Steiner forest}.

For a given three-point set $C=\{A_1, A_2, A_3\}$ a Steiner tree is unique, and,
\begin{itemize}
    \item if the maximal angle of the triangle $\triangle A_1A_2A_3$ is at least $2\pi/3$, i.e. $\angle A_k A_l A_m \geq 2\pi/3$ for some $\{k,l,m\}=\{ 1,2,3\}$, then the Steiner tree is the union of two segments $[A_kA_l] \cup [A_lA_m]$;
    \item otherwise, the Steiner tree is the union of three segments $[A_1F] \cup [A_2F] \cup [A_3F]$, where $F$ is the Fermat--Torricelli point, i.e. the point satisfying
    \[
    \angle A_1FA_2= \angle A_2FA_3=\angle A_3FA_1=2\pi/3.
    \]
    Note that such $F$ belongs to the plane $(A_1A_2A_3)$.
\end{itemize}

A proof of the listed properties of Steiner trees and other interesting information about them one can find in the book~\cite{hwang1992steiner} and in~\cite{gilbert1968steiner}.

\begin{remark}\label{vector}
Let $S$ be a Steiner forest for a finite set $C = \{A_1, \dots, A_m\}$.
For each vertex $A_i$ and for each of its incident edges $X A_i$ of $S$ consider a unit vector codirectional with the vector $\vv{X A_i}$. Then the sum of all these unit vectors is zero. 
\end{remark}

To check Remark~\ref{vector} one can consider a full Steiner tree (as a Steiner forest is a union of full Steiner trees) and prove the statement by induction on the number of Steiner points.

Recall the following widely known facts.

\begin{proposition}\label{3vectors}
There exist at most three unit vectors in $\R^n$ with pairwise angles at least $2\pi/3$. 
In the case of three vectors, they are coplanar, and pairwise angles between them are equal to $2\pi/3$.
Moreover, there exists a constant $c>0$ such that among any four unit vectors in $\R^n$ there exists a pair of vectors with the angle between them at most $2\pi/3-c$.
\end{proposition}
The proof uses the argument from Section 9 of~\cite{gilbert1968steiner}.
\begin{proof}
Let $\vv{v_1}, \dots, \vv{v_k}$ be unit vectors with pairwise angles at least $\alpha$.
Then
\[
0 \leq |\vv{v_1} + \dots + \vv{v_k}|^2 = \sum_i |\vv{v_i}|^2 + \sum_{i \neq j} \langle \vv{v_i}, \vv{v_j} \rangle \leq k + k(k - 1) \cos\alpha,
\]
i.e. $\cos\alpha \geq -\frac{1}{k - 1}$.
Then $\alpha = 2\pi / 3$ implies $k \leq 3$, and in the case $k = 3$ all pairwise angles must be equal to $2\pi / 3$.
Let us further show that in this case the three vectors are coplanar.
We may assume that $\vv{v_1}, \vv{v_2}, \vv{v_3} \in \R^3$ and
\[
\vv{v_1} = (1, 0, 0),\quad \vv{v_2} = (-1/2, \sqrt{3} / 2, 0),\quad \vv{v_3} = (a, b, c).
\]
Then $a = \langle \vv{v_1}, \vv{v_3} \rangle = -1/2$, and $-1/2 = \langle \vv{v_2}, \vv{v_3} \rangle = 1/4 + \sqrt{3}b / 2$, hence $b = -\sqrt{3} / 2$.
It follows that $c^2 = 1 - a^2 - b^2 = 0$, so $\vv{v_1}, \vv{v_2}$ and $\vv{v_3}$ are indeed coplanar.

Finally, if $k = 4$, then $\alpha \leq \arccos(-1/3) = 2\pi / 3 - c$ for some constant $c > 0$.
\end{proof}

\begin{proposition}\label{angles4}
Let $ A, B, C \in \mathbb{R}^n $ be three different points such that $|AB| = |BC| = \varepsilon$ and $\angle ABC \leq 2\pi / 3 - c$ for some $c > 0$. Then 
\[
\H (\St (A, B, C)) \leq (2 - d) \varepsilon,
\]
where $ \St(A, B, C) $ is a Steiner tree for the set $ \{A, B, C \} $ and $ d > 0 $ is a constant depending only on $c$.
\end{proposition}
\begin{proof}
Consider the equilateral triangle $\triangle ACD$ such that $A,B,C,D$ lie in the same plane, and in that plane $B$ and $D$ lie on different sides of the line $(AC)$.
It is known (see, for example, Section 5 in~\cite{gilbert1968steiner}) that $\H(\St(A, B, C)) = |BD|$.
Let $E$ be the intersection of segments $[AC]$ and $[BD]$, then $[BE]$ and $[DE]$ are heights of triangles $\triangle ABC$ and $\triangle ACD$ respectively.
Denote $\alpha := \angle ABC / 2$, then $0 < \alpha \leq \pi / 3 - c / 2$, and
\[
|BD| = |BE| + |DE| = \varepsilon\cos\alpha + \sqrt{3}\varepsilon\sin\alpha = 2\varepsilon \sin(\alpha + \pi / 6)\leq 2\varepsilon \cos(c/2) \leq (2 - d) \varepsilon
\]
for some $d > 0$ which depends only on $c$.
\end{proof}


\begin{remark}
There exists a constant $c_1>0$ such that for any $X \in \R^n$, $\rho > 0$ and $m\geq 4$ different points $A_1, \dots, A_m$ on the sphere $\partial B_{\rho}(X)$ a Steiner tree for the set $\{X, A_1, \dots, A_m\}$ has the length not greater than $(m-c_1)\rho$.
\label{no_inf}
\end{remark}
\begin{proof}
Due to Proposition~\ref{3vectors} we may assume that $\angle A_1XA_2 \leq 2\pi / 3 - c$ for some $c > 0$.
Then, by Proposition~\ref{angles4}, the length of a Steiner tree $\St(X, A_1, A_2)$ for the set $\{X, A_1, A_2\}$ is at most $(2 - d) \rho$ for some $d > 0$ which depends only on $c$, i.e. does not depend on the choice of points $A_1, \dots, A_m$.
Finally, the connected set $\St(X, A_1, A_2) \cup [XA_3] \cup \dots \cup [XA_m]$ has the length at most $(m - d)\rho$, and a Steiner tree for the set $\{X, A_1, \dots, A_m\}$ is not longer than that set.
\end{proof}



\section{Notation and definitions}\label{not}

\begin{itemize}
\item For a given set $X \subset \mathbb{R}^n$ the sets $\overline{X}$, $\Int (X)$ and $\partial X$ represent its closure, interior and boundary, respectively, and $\conv X$ denotes its convex hull.

\item Let $B, C \in \mathbb{R}^n$ be points.
Then $|BC|$ denotes Euclidean distance between $B$ and $C$.
Also, $[BC]$ denotes the segment with endpoints at $B$ and $C$, $\vv{BC}$ denotes the vector starting at $B$ and ending at $C$, both $[BC)$ and $(CB]$ denote the ray which starts at $B$ and contains $C$, and $(BC)$ denotes the line which contains $B$ and $C$.

\item The symbol~$\angle$ is applicable to two lines, two rays, two vectors or three points.
It denotes the unique angle between the objects (in case of three points, $\angle ABC = \angle ( \vv{BA}, \vv{BC} )$) belonging to $[0, \pi]$  for rays, vectors and points, and to $[0, \pi / 2]$ for lines.

\item Given a sequence of rays $(a_kb_k]$, we say that $(a_kb_k] \rightarrow (ab]$ if $b_k \rightarrow b$ and $\angle((a_kb_k], (ab]) \rightarrow 0$.

\item We will say that two vectors, rays or lines are \textit{$\gamma$-orthogonal} to each other if the angle between them differs from $\pi / 2$ by at most $\gamma$.
Similarly, we will say that two vectors, rays or lines are \textit{$\gamma$-parallel} to each other if the angle between them differs from $0$ by at most $\gamma$.
We will sometimes abuse the notation and say that a ray $[AB)$ or a vector $\vv{AB}$ is $\gamma$-orthogonal ($\gamma$-parallel) to a line $(CD)$ if the line $(AB)$ is $\gamma$-orthogonal ($\gamma$-parallel) to $(CD)$.

\item Given two vectors $\vv{a}, \vv{b} \in \mathbb{R}^n$, the dot product between them is denoted by $\langle \vv{a}, \vv{b} \rangle \defeq \sum_{i = 1}^n a_i b_i$.


\item The supremum of the distance between a pair of points from a set $X \subset \mathbb{R}^n$ will be referred to as the \textit{diameter} of $X$:
\[
\diam X \defeq \sup_{a, b \in X} \dist(a, b).
\]

\item For a (usually finite) set $X \subset \mathbb{R}^n$, its cardinality, i.e. the number of elements in the set, is denoted by $ \sharp X $.

\item We will write $f(\varepsilon)=o_\varepsilon(g(\varepsilon))$ if and only if $\frac{f(\varepsilon)}{g(\varepsilon)}\rightarrow 0$ as $\varepsilon \rightarrow 0$. In particular, $f(\varepsilon)=o_\varepsilon(1)$ if and only if $f(\varepsilon) \rightarrow 0$ as $\varepsilon \rightarrow 0$. If the dependence on $\varepsilon$ is clear from the context, we will often write $o(\cdot)$ instead of $o_\varepsilon(\cdot)$.
\end{itemize}


\begin{remark}\label{ord_pos}
If $\Gamma \subset \mathbb{R}^n$ is a connected set of a positive length, then $ \oord_x \Gamma > 0 $ if and only if $ x \in \overline{\Gamma} $.
\end{remark}
\begin{proof}
If $x \not\in \overline{\Gamma}$, then $\dist(x, \Gamma) > 0$, hence $\sharp (\Gamma \cap \partial U) = 0$ for any open $U \ni x$ with $\diam U < \dist(x, \Gamma)$, so $\oord_x \Gamma = 0$.
If $x \in \overline{\Gamma}$, then there exists $y \in \Gamma$, $y \neq x$.
Then $\sharp (\Gamma \cap \partial U) > 0$ for any open $U \ni x$ with $\diam U < \dist(x, y)$, because otherwise two open (in the relative topology of $\Gamma$) sets $\Gamma \cap U$ and $\Gamma \setminus \overline{U}$ would form a partition of $\Gamma$, which is not possible since $\Gamma$ is connected.
\end{proof}

\begin{remark}\label{ord_comp}
Let $\Gamma \subset \mathbb{R}^n$ be a closed connected set and let the order of a point $x \in \Gamma$ be finite: ${\oord_x \Gamma < \infty}$.
By Remark~\ref{ord_pos}, $\oord_x \Gamma_1 > 0$ for any connected component $\Gamma_1$ of $\Gamma \setminus \{x\}$, so $\oord_x \Gamma$ is at least the number of such connected components.
Thus, the number of connected components of $\Gamma \setminus \{x\}$ is finite for any $x \in \Gamma$ with $\oord_x \Gamma < \infty$.
\end{remark}


\begin{definition}
We define the circular order $\oordball_x \Gamma$ of a point $x \in \R^n$ in a set $\Gamma \subset \mathbb{R}^n$ as the supremum of $m \in \N$ for which there exists $ \rho = \rho(x) > 0 $ such that for any positive $\varepsilon < \rho$ the inequality $ \sharp (\Gamma \cap \partial B_\varepsilon (x)) \geq m $ holds.
In other words,
\[
\oordball_x \Gamma := \sup\{m \in \mathbb{N}\ |\ \exists \rho>0: \sharp(\Gamma \cap \partial B_\varepsilon(x))\geq m\ \forall \text{ positive } \varepsilon < \rho \}.
\]
\end{definition}

\begin{remark}\label{rm:ordball}
Let the circular order of a point $x \in \R^n$ in a set $\Gamma \subset \mathbb{R}^n$ be finite: $m \defeq \oordball_x \Gamma < \infty$. Then there exists $ \rho = \rho(x) > 0 $ such that for any positive $\varepsilon < \rho$ the inequality $ \sharp (\Gamma \cap \partial B_\varepsilon (x)) \geq m $ holds, and for any positive $\varepsilon < \rho $ there exists a positive $\delta < \varepsilon $ such that $ \sharp (\Gamma \cap \partial B_{\delta} (x)) = m $.
\end{remark}

\begin{remark}\label{rem_ord}
For any set $\Gamma \subset \mathbb{R}^n$ and any $x \in \R^n$ the inequality $\oordball_x \Gamma \geq \oord_x \Gamma$ holds.
\end{remark}
\begin{proof}
The supremum in the definition of $\oordball$ is on a (not strictly) larger family than in the definition of $\oord$.
\end{proof}

\begin{proposition}\label{coarea}
Let $\Gamma \subset \R^n$ be a set with $\H(\Gamma) < \infty$, let $x \in \R^n$, $k \geq 0$, $0 \leq \varepsilon_1 < \varepsilon_2$, and let $\sharp(\Gamma \cap \partial B_\varepsilon(x)) \geq k$ for any $\varepsilon_1 < \varepsilon < \varepsilon_2$. Then
\[
\H( \Gamma \cap \overline{B_{\varepsilon_2}(x) \setminus B_{\varepsilon_1}(x)} ) \geq k (\varepsilon_2 - \varepsilon_1).
\]
\end{proposition}
\begin{proof}
We apply Theorem 1.1 from~\cite{esmayli2021coarea}, i.e. the coarea inequality, also known as Eilenberg's inequality~\cite{eilenberg1938phi}.
Substituting $s = t = 1$, $X = \mathbb R^n$, $Y = \R$, $f(\cdot) = \dist(x, \cdot)$, $E = \Gamma \cap \overline{B_{\varepsilon_2}(x) \setminus B_{\varepsilon_1}(x)}$ into the mentioned theorem, we get
\[
\H(\Gamma \cap \overline{B_{\varepsilon_2}(x) \setminus B_{\varepsilon_1}(x)}) \geq \int_{\varepsilon_1}^{\varepsilon_2} \HH^0(\Gamma \cap \partial B_\varepsilon(x)) d\H(\varepsilon) = \int_{\varepsilon_1}^{\varepsilon_2} \sharp(\Gamma \cap \partial B_{\varepsilon}(x)) d\H(\varepsilon) \geq k(\varepsilon_2 - \varepsilon_1).
\]
\end{proof}

\begin{remark}\label{geq_ord}
Let $\Gamma \subset \mathbb{R}^n$ be a set with $\H(\Gamma) < \infty$ and let the circular order of a point $x\in \R^n$ be finite: ${\oordball_x \Gamma < \infty}$. Then there exists $\rho = \rho(x) > 0$ such that for any positive $\varepsilon < \rho$ the inequality 
\[
\H(\Gamma \cap B_\varepsilon(x))\geq \oordball_x\Gamma \cdot \varepsilon
\]
holds.
\begin{proof}
Apply Proposition~\ref{coarea} with $k = \oordball_x \Gamma$, $\varepsilon_1 = 0$ and $\varepsilon_2 = \rho$, where $\rho$ is from Remark~\ref{rm:ordball}.
\end{proof}
\end{remark}


\begin{definition}
Let $\Sigma$ be a local minimizer for a compact set $M \subset \mathbb{R}^n$ and $r > 0$.
A point $x \in \Sigma$ is energetic if for every $\varepsilon > 0$ the inequality
\[
\F_{M}(\Sigma \setminus B_{\varepsilon}(x)) > \F_{M}(\Sigma)
\]
holds.
\end{definition}

The set of all energetic points of $\Sigma$ is denoted by $\G_\Sigma$.
Each local minimizer $\Sigma$ can be split into three disjoint subsets:
\[
\Sigma=\E_\Sigma \sqcup \X_\Sigma \sqcup \S_\Sigma,
\]
where $\X_\Sigma \subset \G_\Sigma$ is the set of \textit{isolated energetic points} (i.e. every $x\in \X_\Sigma$ is energetic and there exists $\rho = \rho(x) > 0$ such that $\G_\Sigma \cap B_\rho(x) = \{x\}$), $\E_\Sigma \defeq \G_\Sigma\setminus \X_\Sigma$ is the set of \textit{non-isolated energetic points} and $\S_\Sigma \defeq \Sigma \setminus \G_\Sigma$ is the set of \textit{non-energetic points}.

For every energetic point $x \in \G_\Sigma$ there exists a \textit{point, corresponding to $x$}, i.e. a point $y(x) \in M$ such that $\dist (x, y(x)) = r$ and $\Sigma \cap B_r(y(x)) = \emptyset$.
We will refer to this property as \textit{the basic property of energetic points}.
Note that an energetic point $x \in \G_\Sigma$ may have several points corresponding to it, and a non-energetic point $x \in \S_\Sigma$ may also have a corresponding point.

\begin{remark}\label{rem:non-energ}
Let $\Sigma$ be a local minimizer for a compact set $M \subset \mathbb{R}^n$ and $r > 0$, and let $x \in \S_\Sigma$ be a non-energetic point.
Then $\oord_x \Sigma = \oordball_x \Sigma \in \{2, 3\}$, and there exists $\rho = \rho(x) > 0$ such that
\begin{itemize}
\item if $\oord_x \Sigma = 2$, then $\Sigma \cap \overline{B_\rho(x)}$ is a segment with endpoints at $\partial B_\rho(x)$ which contains $x$;
\item if $\oord_x \Sigma = 3$, then $\Sigma \cap \overline{B_\rho(x)}$ is a regular tripod centered at $x$, i.e. the union of three segments with one endpoint at $x$ and another at $\partial B_\rho(x)$, and pairwise angles $2 \pi / 3$ between them.
\end{itemize}
\end{remark}
\begin{proof}
Since the length of $\Sigma$ is finite, there exists an arbitrarily small $\varepsilon > 0$ such that $\sharp (\Sigma \cap \partial B_\varepsilon(x)) < \infty$.
Let this $\varepsilon$ be small enough that $\F_{M}(\Sigma \setminus B_\varepsilon(x)) = \F_{M}(\Sigma)$ (due to the fact that $x$ is non-energetic) and that $\Sigma$ is the shortest set among all sets $\Sigma'$ with $\diam(\Sigma \triangle \Sigma') \leq \varepsilon$ and $\F_M(\Sigma') \leq r$ (due to the definition of local minimizers).
Thus the set $\Sigma \cap \overline{B_\varepsilon(x)}$ must be a Steiner forest for a finite set of points $\Sigma \cap \partial B_\varepsilon(x)$.
Then the claim of the remark is true for some $\rho \leq \varepsilon$.
\end{proof}

We will call a non-energetic point $x \in \S_\Sigma$ with $\oord_x \Sigma = 3$ a \textit{branching point} or a \textit{Steiner point}.
We will also call a branching point any energetic point $x \in \G_\Sigma$ with $\oord_x \Sigma > 2$ (later we will show that in this case $\oord_x \Sigma = 3$, and, although $x$ is not necessarily the center of a regular tripod in $\Sigma$, tangent rays of $\Sigma$ at $x$ do form a regular tripod).


Note that it is possible for a (local) minimizer to have no non-energetic points at all.
Moreover, in some sense, any (local) minimizer does not have any non-energetic points:

\begin{example}\label{rem:all-energ}
Let $\Sigma$ be a (local) minimizer for a compact set $M \subset \mathbb{R}^n$ and $r > 0$.
Then $\hat \Sigma \defeq \Sigma \times \{0\} \subset \mathbb{R}^{n + 1}$ is a (local) minimizer for $\hat M = (M \times \{0\}) \cup (\Sigma \times \{r\}) \subset \mathbb{R}^{n + 1}$ and $\E_{\hat \Sigma} = \hat \Sigma$.
\end{example}
\begin{proof}
We consider the case of minimizers; the same reasoning works in the case of local minimizers.
By definition, $\F_{\hat M}(\hat \Sigma) \leq r$.
Suppose $\hat \Sigma$ is not a minimizer, then there exists a connected set $\Gamma \subset \R^{n + 1}$ with $\F_{\hat M}(\Gamma) \leq r$ and $\H(\Gamma) < \H(\hat \Sigma)$.
Let $\Gamma'$ be the orthogonal projection of $\Gamma$ onto $\R^n \times \{0\}$, then $\H(\Gamma') \leq \H(\Gamma) < \H(\hat \Sigma) = \H(\Sigma)$ and $\F_{M \times \{0\}}(\Gamma') \leq r$, which contradicts $\Sigma$ being a minimizer for $M$.
\end{proof}

\begin{remark}\label{emptyball}
Let $\Sigma$ be a local minimizer for a compact set $M \subset \mathbb{R}^n$ and $r > 0$, and let $x_k \in \G_\Sigma$ be a sequence of energetic points converging to a non-isolated energetic point $x\in \E_\Sigma$: $x_k \rightarrow x$.
By the basic property of energetic points, for every point $x_k \in \G_\Sigma$ there exists a corresponding point $y_k = y_k(x_k) \in M$. Let $y$ be an arbitrary limit point of the set $\{y_k\}$. Then $y \in M$ and $\Sigma \cap B_r(y) = \emptyset$, i.e. $y$ corresponds to $x$.
\end{remark}

\begin{definition}\label{def:tgray}
We will say that the ray $ (ax] $ is a \textit{tangent ray} of a set $ \Gamma \subset \mathbb{R}^n $
at a point $ x\in \Gamma $ if there exists a
sequence of points $ x_k \in \Gamma \setminus \{x\}$ such that $ x_k \rightarrow x $ and $ \angle x_kxa \rightarrow 0 $.
\end{definition}

\begin{definition}\label{def:one-sided}
We will say that the ray $ (ax] $ is a \textit{one-sided tangent} of a set $ \Gamma  \subset \mathbb{R}^n $ at a point $ x \in \Gamma $ if there exists a connected component $\Gamma_1$ of $\Gamma \setminus \{x\}$ such that $x \in \overline{\Gamma_1}$ and that any sequence of points $x_k \in \Gamma_1$ with the property $x_k \rightarrow x$ satisfies $\angle x_kxa \rightarrow 0$.
In this case we will also say that $(ax]$ is tangent to the connected component $\Gamma_1$.
\end{definition}

\begin{remark}
If $(ax]$ is a one-sided tangent of a set $\Gamma \subset \mathbb{R}^n$ at a point $x \in \Gamma$, then it is also a tangent ray of $\Gamma$ at $x$, but not vice versa.
\end{remark}

 \section{Results} \label{res}

\begin{definition}
Let $\Sigma$ be a local minimizer for a compact set $M \subset \mathbb{R}^n$ and $r > 0$.
Given any $x \in \Sigma$ and any $ \varepsilon > 0 $, define the set
\[
U(x,\varepsilon) \defeq \left \{ u\in M\ \big |\ \dist(u, \Sigma \cap \overline{B_{\varepsilon}(x)})\leq r,\ \dist(u, \Sigma \setminus \overline{B_{\varepsilon}(x)}) \geq r \right \}.
\]
\end{definition}

\begin{remark}\label{rem_c}
Let $\Sigma$ be a local minimizer for a compact set $M \subset \mathbb{R}^n$ and $r > 0$.
For any $x \in \Sigma$ the set $\bigcap_{\varepsilon > 0} U(x, \varepsilon)$ is precisely the set of corresponding points. In particular, if $U(x,\varepsilon)=\emptyset$ for some $\varepsilon > 0$, then a point $x \in \Sigma$ is non-energetic.
\end{remark}

\begin{remark}\label{rem_cC}
Let $\Sigma$ be a local minimizer for a compact set $M \subset \mathbb{R}^n$ and $r > 0$, and let $x \in \G_\Sigma$ be an energetic point.
Then for any $u \in U(x, \varepsilon)$ there exists a corresponding point $y(x)$ such that $\angle uxy(x) = o_\varepsilon(1)$.
\end{remark}
\begin{proof}
Consider the contrary, then there exists a constant $c > 0$ and sequences $\varepsilon_i \rightarrow 0$ and $u_i \in U(x, \varepsilon_i)$ such that $\angle u_i xy(x) \geq c$ for any corresponding point $y(x)$.
Due to compactness, there exists a converging subsequence $u_{i_k} \rightarrow y$ for some limit point $y$.
By Remark~\ref{rem_c}, $y \in \bigcap_{\varepsilon > 0} U(x, \varepsilon)$ corresponds to $x$.
But then $\angle u_{i_k} xy \geq c$ for any $k$, so $0 = \angle yxy \geq c$, which is a contradiction.
\end{proof}

\definecolor{rvwvcq}{rgb}{0.08235294117647059,0.396078431372549,0.7529411764705882}
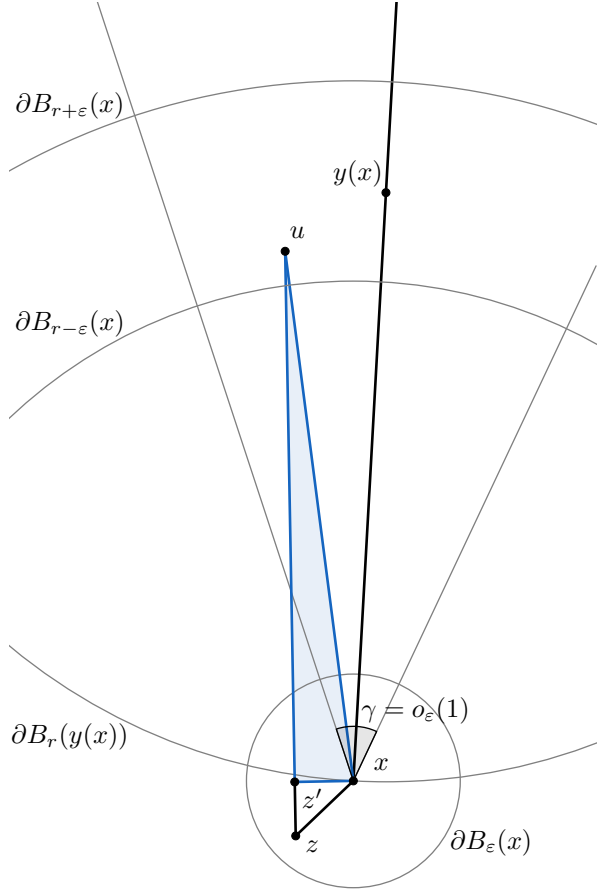
\begin{figure}[ht]
\center{}
\begin{tikzpicture}[line cap=round,line join=round,>=triangle 45,x=1.0cm,y=1.0cm]
\clip(-6.2,-8.17) rectangle (1.5934,4.3656);
\fill[line width=1.pt,color=rvwvcq,fill=rvwvcq,fill opacity=0.10000000149011612] (-1.6334,-5.9862) -- (-2.412555781431965,-6.0002644730924475) -- (-2.54,1.06) -- cycle;
\draw [line width=0.5pt,color=gray] (-1.6334,-5.9862) circle (1.4211016853131937cm);
\draw [line width=0.5pt,color=gray] (-1.2,1.84) circle (7.838191245434115cm);
\draw [line width=1.pt,domain=-1.6334:1.5934000000000035] plot(\x,{(--10.188896--7.8262*\x)/0.4334});
\draw [line width=0.5pt,domain=-1.6334:1.5934000000000035,color=gray] plot(\x,{(-7.738091999999999--6.5262*\x)/3.0734});
\draw [line width=0.5pt,domain=-5.6423999999999985:-1.6334,color=gray] plot(\x,{(--28.500647999999998--7.92*\x)/-2.6});
\draw [line width=1.pt] (-2.412555781431965,-6.0002644730924475)-- (-2.3996,-6.718);
\draw [line width=1.pt,color=rvwvcq] (-1.6334,-5.9862)-- (-2.412555781431965,-6.0002644730924475);
\draw [line width=1.pt,color=rvwvcq] (-2.412555781431965,-6.0002644730924475)-- (-2.54,1.06);
\draw [line width=1.pt,color=rvwvcq] (-2.54,1.06)-- (-1.6334,-5.9862);
\draw [line width=1.pt,color=black] (-2.3996,-6.718)-- (-1.6334,-5.9862);
\draw [line width=0.5pt,color=gray] (-1.6334,-5.9862) circle (6.649211832991937cm);
\draw [line width=0.5pt,color=gray] (-1.6334,-5.9862) circle (9.314100547020093cm);
\draw [fill=black] (-1.6334,-5.9862) circle (1.5pt);
\draw[color=black] (-1.2622,-5.7685) node {$x$};
\draw[color=black] (0.2,-6.8) node {$\partial B_\varepsilon(x)$};
\draw[color=black] (-0.80598,-5.05) node {$\gamma = o_\varepsilon(1)$};
\draw [fill=black] (-2.54,1.06) circle (1.5pt);
\draw[color=black] (-2.3754,1.3079) node {$u$};
\draw [fill=black] (-2.412555781431965,-6.0002644730924475) circle (1.5pt);
\draw[color=black] (-2.1818,-6.2203) node {$z'$};
\draw [fill=black] (-2.3996,-6.718) circle (1.5pt);
\draw[color=black] (-2.1818,-6.8508443) node {$z$};
\draw[color=black] (-5.4,0.1) node {$\partial B_{r-\varepsilon}(x)$};
\draw [fill=black] (-1.2,1.84) circle (1.5pt);
\draw[color=black] (-1.572,2.0965) node {$y(x)$};
\draw[color=black] (-5.4,3) node {$\partial B_{r+\varepsilon}(x)$};
\draw[color=black] (-5.4,-5.4) node {$\partial B_r(y(x))$};
\draw [shift={(-1.6334,-5.9862)},line width=0.5pt,color=black,fill=black,fill opacity=0.10000000149011612] (0,0) -- (64.78272748171172:0.726) arc (64.78272748171172:108.17412006546908:0.726) -- cycle;
\end{tikzpicture}
\caption{Illustration to the proof of Lemma~\ref{Ooo}}
\label{l318}
\end{figure}

Note that for any $u\in U(x, \varepsilon)$ the inequality $|x u|\leq r+\varepsilon$ holds.
In the following lemma we show that $\varepsilon$ on the right hand side can be replaced by $o(\varepsilon)$.

\begin{lemma}\label{Ooo}
Let $\Sigma$ be a local minimizer for a compact set $M \subset \mathbb{R}^n$ and $r > 0$, and let $x \in \G_\Sigma$ be an energetic point.
Then for any $u\in U(x, \varepsilon)$ the inequality $|xu|\leq r+o(\varepsilon)$ holds.
\end{lemma}
\begin{proof}
Consider an arbitrary point $u\in U(x, \varepsilon)$ and let a point $z\in \Sigma \cap \overline{B_\varepsilon(x)}$ be such that $|uz|\leq r$ (see Figure~\ref{l318}).
Let $y(x)$ be a corresponding point such that $\angle ux y(x) = o_\varepsilon(1)$ as in Remark~\ref{rem_cC}.
By the definition of a corresponding point, $z\notin B_r(y(x))$. 
Then
\[
\angle zxu \geq \angle zxy(x) - \angle uxy(x) \geq \pi/2 - o_\varepsilon(1), \quad\text{so}\quad
\angle uzx = \pi - \angle zxu - \angle zux\leq \pi/2+o_\varepsilon(1).
\]
Define $z'$ as the base of the perpendicular from $x$ to $(zu)$, if $\angle uzx<\pi/2$, and as equal to $z$, otherwise. Then $\angle uz'x = \pi/2 + o_\varepsilon(1)$.

Consider the triangle $\triangle z'xu$. In this triangle $|z'u|\leq |zu|\leq r$, $|xz'|\leq |xz| \leq \varepsilon$ and $\angle uz'x = \pi/2 + o_\varepsilon(1)$. 
Then, by the law of cosines,
\[
|xu|=\sqrt{|xz'|^2+|z'u|^2 - 2\cos \angle uz'x \cdot |z'u| \cdot |xz'|}\leq \sqrt{r^2+\varepsilon^2+\varepsilon\cdot o_\varepsilon(1)}=r+o(\varepsilon).
\]
\end{proof}

\begin{lemma}\label{cube}
Let $r > 0$. There exists $c = c(n) > 0$ and $\rho = \rho(n, r) > 0$ such that for any positive $\varepsilon < \rho$ and any $x \in \mathbb{R}^n$ there exists a connected set $S(x, \varepsilon) \ni x$ satisfying $\overline{B_{r + \varepsilon}(x)} \subset \overline {B_r(S(x, \varepsilon))}$ and $\H(S(x, \varepsilon)) = c\varepsilon$.
\end{lemma}
\begin{proof}
Let $x = (x_1, \dots, x_n)$, let $c' \defeq 4\sqrt{n}$, and let $V$ be the set consisting of $2^n$ points which differ from $x$ by $c'\varepsilon$ at each coordinate.
Let $S = S(x, \varepsilon)$ be a Steiner tree for the set $V \cup \{x\}$.
Note that $S$ is connected, $x \in S$ and $\H(S) = c \varepsilon$ for some $c = c(n)$.

Let $y = (y_1, \dots, y_n) \in \overline{B_{r + \varepsilon}(x)}$.
Due to symmetry, we can assume that $l_i \defeq y_i - x_i \geq 0$ for $1 \leq i \leq n$.
Let $l \defeq \dist(y, x)$ and note that $0 \leq l \leq r + \varepsilon$ and that $l_i \geq l / \sqrt{n}$ for some $i$.
Let $v \defeq (x_1 + c'\varepsilon, \dots, x_n + c'\varepsilon)$.
Note that $v \in V$, so $v \in S$.
Let us estimate the distance from $y$ to $v$:
\[
(\dist(y, v))^2 = \sum_{i = 1}^n (l_i - c'\varepsilon)^2 = l^2 - 2c'\varepsilon \sum_{i = 1}^n l_i + o(\varepsilon) \leq l^2 - \frac{2lc'}{\sqrt{n}}\varepsilon + o(\varepsilon) = l^2 - 8l\varepsilon + o(\varepsilon).
\]
If $l < \frac{r}{2}$, then $(\dist(y, v))^2 \leq \frac{r^2}{4} + o(\varepsilon) \leq r^2$ for a sufficiently small $\varepsilon$.
On the other hand, if $l \geq \frac{r}{2}$, then
\[
(\dist(y, v))^2 \leq (r + \varepsilon)^2 - 4r\varepsilon + o(\varepsilon) = r^2 - 2r\varepsilon + o(\varepsilon) \leq r^2
\]
if $\varepsilon$ is small enough.
This concludes the proof of the lemma.
\end{proof}

By combining Lemma~\ref{Ooo} and Lemma~\ref{cube}, we get the following corollary.

\begin{corollary}\label{main_c}
Let $\Sigma$ be a local minimizer for a compact set $M \subset \mathbb{R}^n$ and $r > 0$, and let $x \in \Sigma$.
Then there exists $\rho = \rho(x) > 0$ such that for any positive $\varepsilon < \rho$ there exists a connected set $S^*(x, \varepsilon) \ni x$ satisfying $\H(S^*(x, \varepsilon))=o(\varepsilon)$ and
\[
\F_{M}(\Sigma \setminus B_\varepsilon( x) \cup S^*(x, \varepsilon)) \leq \F_M(\Sigma).
\]
\end{corollary}

\begin{corollary}
Let $\Sigma$ be a local minimizer for a compact set $M \subset \mathbb{R}^n$ and $r > 0$.
Then $\oordball_x\Sigma \leq 3$ for any $x \in \R^n$.
\end{corollary}
\begin{proof}
Assume the contrary, i.e. $\oordball_x \Sigma \geq 4$.
Since the length of $\Sigma$ is finite, there exists an arbitrarily small $\varepsilon > 0$ such that $\sharp ( \Sigma \cap \partial B_\varepsilon(x) ) < \infty$.
Using this fact if $\oordball_x \Sigma = \infty$, and Remark~\ref{rm:ordball} otherwise, we see that there exists a sequence of positive $\varepsilon_k \rightarrow 0$ such that $m_k \defeq \sharp(\Sigma \cap \partial B_{\varepsilon_k}(x)) \geq 4$ is finite and $\sharp(\Sigma \cap \partial B_{\delta}(x)) \geq m_k$ for every positive $\delta\leq\varepsilon_k$.
Due to Proposition~\ref{coarea}, $\H(\Sigma \cap B_{\varepsilon_k}(x)) \geq m_k \cdot \varepsilon_k$.
On the other hand, if $S_k$ is a Steiner tree for the set of points $(\Sigma \cap \partial B_{\varepsilon_k}(x)) \cup \{x\}$, then, by Remark~\ref{no_inf}, $\H(S_k) \leq (m_k - c_1) \cdot \varepsilon_k$ for some constant $c_1 > 0$.
Then, for a sufficiently big $k$, the set 
\[
\Sigma \setminus B_{\varepsilon_k}(x) \cup S^*(x, \varepsilon_k) \cup S_k,
\]
where $S^*(x, \varepsilon_k)$ is from Corollary~\ref{main_c}, is connected, has the energy not greater than $\Sigma$ and is shorter than $\Sigma$, which is not possible.
\end{proof}

\begin{definition}
Let $\Gamma \subset \mathbb{R}^n$ be a closed connected set and let the circular order of a point $x \in \Gamma$ be finite: $\oordball_x \Gamma < \infty$.
We will say that $\varepsilon > 0$ is $(\Gamma, x)$-nice if $\sharp (\Gamma \cap \partial B_{\varepsilon} (x)) = \oordball_x \Gamma$ and for every connected component $\Gamma_1$ of $\Gamma \setminus \{x\}$ the inequality $\sharp (\Gamma_1 \cap \partial B_{\varepsilon}(x)) > 0$ holds.
\end{definition}

Note that, in view of Remark~\ref{rm:ordball}, for any closed connected set $\Gamma \in \R^n$ and any $x \in \Gamma$ with $\oordball_x \Gamma < \infty$ there exists a sequence of $(\Gamma, x)$-nice $\varepsilon_k$ such that $\varepsilon_k \rightarrow 0$.

\begin{lemma}\label{lm:tech}
Let $\Sigma$ be a local minimizer for a compact set $M \subset \mathbb{R}^n$ and $r > 0$, and let $x \in \Sigma$.
There exists $\rho = \rho(x) > 0$ such that
\begin{enumerate}
\item \label{tech:len} for any positive $\varepsilon < \rho$ the equality $\H(\Sigma \cap B_\varepsilon(x)) = \oordball_x\Sigma \cdot \varepsilon + o(\varepsilon)$ holds;
\item \label{tech:angle} for any $(\Sigma, x)$-nice $\varepsilon < \rho$ and any $A, B \in \Sigma \cap \partial B_\varepsilon(x)$, $A \neq B$, the inequality $\angle AxB \geq 2 \pi / 3 - o_\varepsilon(1)$ holds;
\item \label{tech:nice} for any $(\Sigma, x)$-nice $\varepsilon < \rho$ there exists $(\Sigma, x)$-nice $\delta < \varepsilon$ such that $\delta = \varepsilon \cdot (1 - o_\varepsilon(1))$.
\end{enumerate}
\end{lemma}
\begin{proof}\leavevmode
\begin{enumerate}
\item Let us start by showing that the claim holds for any $(\Sigma, x)$-nice $\varepsilon < \rho$.
Indeed, assume the contrary, then, in view of Remark~\ref{geq_ord}, there exists $c > 0$ and a sequence of $(\Sigma, x)$-nice $\varepsilon_k > 0$ such that $\varepsilon_k \rightarrow 0$ and $\H( \Sigma \cap B_{\varepsilon_k}(x) ) \geq (\oordball_x\Sigma + c) \cdot \varepsilon_k$ for every $k$.
Then, for a sufficiently big $k$, the set 
\[
\Sigma \setminus B_{\varepsilon_k}(x) \cup S^*(x, \varepsilon_k) \cup \bigcup_{A\in \Sigma \cap \partial B_{\varepsilon_k}(x)} [Ax],
\]
where $S^*(x, \varepsilon_k)$ is from Corollary~\ref{main_c}, is connected, has the energy not greater than $\Sigma$ and is shorter than $\Sigma$, which is not possible.

Now, assume that the claim is false if $\varepsilon < \rho$ is allowed not to be $(\Sigma, x)$-nice.
Then there exists $c > 0$ and a sequence $\varepsilon_k > 0$ such that $\varepsilon_k \rightarrow 0$, $\H( \Sigma \cap B_{\varepsilon_k}(x) ) \geq (\oordball_x \Sigma + c) \cdot \varepsilon_k$ and $\sharp (\Sigma \cap \partial B_{\varepsilon_k}(x)) \geq \oordball_x \Sigma + 1$.
Let $\delta_k \geq \varepsilon_k$ be the infimum of $(\Sigma, x)$-nice $\delta > \varepsilon_k$.
On the one hand, $\delta_k$ is a limit point of some sequence of $(\Sigma, x)$-nice $\varepsilon < \rho$, so $\H(\Sigma \cap B_{\delta_k}(x)) = \oordball_x\Sigma \cdot \delta_k + o(\delta_k)$.
On the other hand, due to Proposition~\ref{coarea},
\[
\H(\Sigma \cap B_{\delta_k}(x)) \geq (\oordball_x\Sigma + c) \cdot \varepsilon_k + ( \oordball_x\Sigma + 1 ) \cdot (\delta_k - \varepsilon_k) \geq (\oordball_x\Sigma + \min(c, 1)) \cdot \delta_k
\]
for every $k$.
We have a contradiction.

\item Assume the contrary, i.e. that there exists $c > 0$, a sequence of $(\Sigma, x)$-nice $\varepsilon_k$, and $A_k, B_k \in \Sigma \cap \partial B_{\varepsilon_k}(x)$, $A_k \neq B_k$ for each $k$, such that $\varepsilon_k \rightarrow 0$ and $\angle A_kxB_k < 2 \pi / 3 - c$.
Then, for a sufficiently big $k$, the set
\[
\Sigma \setminus B_{\varepsilon_k}(x) \cup S^*(x,\varepsilon_k) \cup \St(A_k, x, B_k) \cup \bigcup_{C \in ( \Sigma \cap \partial B_{\varepsilon_k}(x) ) \setminus \{ A_k, B_k \}} [Cx],
\]
where $S^*(x,\varepsilon_k)$ is from Corollary~\ref{main_c} and $\St(A_k, x, B_k)$ is a Steiner tree for the set $\{ A_k, x, B_k\}$, is connected and has the energy not greater than $\Sigma$.
In view of Proposition~\ref{angles4} and item~\ref{tech:len}, it is also shorter than $\Sigma$ if $k$ is sufficiently big.
It is not possible, so our assumption is false.

\item Assume the contrary, then there exists $c > 0$ and a sequence of $(\Sigma, x)$-nice $\varepsilon_k$ such that $\varepsilon_k \rightarrow 0$ and $\sharp ( \Sigma \cap \partial B_\delta(x) ) \geq \oordball_x\Sigma + 1$ for every $\delta$ satisfying $(1 - c) \cdot \varepsilon_k \leq \delta < \varepsilon_k$.
Let $0 < \delta_k \leq (1 - c) \cdot \varepsilon_k$ be the supremum of $(\Sigma, x)$-nice $\delta < \varepsilon_k$.
Then, due to Proposition~\ref{coarea},
\[
\H( \Sigma \cap B_{\varepsilon_k}(x) ) \geq \oordball_x\Sigma \cdot \delta_k + ( \oordball_x \Sigma + 1 ) \cdot (\varepsilon_k - \delta_k) \geq ( \oordball_x \Sigma + c ) \cdot \varepsilon_k
\]
for every $k$, which contradicts item~\ref{tech:len}.
\end{enumerate}
\end{proof}

\begin{corollary}\label{cor:struct}
Let $\Sigma$ be a local minimizer for a compact set $M \subset \mathbb{R}^n$ and $r > 0$, and let $x \in \Sigma$.
\begin{enumerate}
\item \label{struct:angle} The angle between any pair of one-sided tangents of $\Sigma$ at $x$ is at least $2 \pi / 3$.
In view of Proposition~\ref{3vectors}, all one-sided tangents of $\Sigma$ at $x$ are coplanar.
\item \label{struct:ord} The equality $\oord_x \Sigma = \oordball_x \Sigma$ holds.
Moreover, $1 \leq \oord_x \Sigma \leq 3$.
\item \label{struct:tang-comp} There are exactly $\oord_x \Sigma$ one-sided tangents of $\Sigma$ at $x$ and exactly $\oord_x \Sigma$ connected components of $\Sigma \setminus \{x\}$ with one-to-one correspondence between them, i.e. every connected component has a one-sided tangent, and every one-sided tangent is tangent to exactly one connected component.
It follows that the same is true if we replace ``one-sided tangents'' with ``tangent rays''.
\item \label{struct:comp} Let $\Sigma_1$ be an arbitrary connected component of $\Sigma \setminus \{x\}$.
Then $\oord_x \Sigma_1 = \oordball_x \Sigma_1 = 1$, and for any $\varepsilon > 0$ the equality $\H(\Sigma_1 \cap B_\varepsilon(x)) = \varepsilon + o(\varepsilon)$ holds.
\item \label{struct:loops} $\Sigma$ does not contain loops (homeomorphic images of circles).
\end{enumerate}
\end{corollary}
\begin{proof}
Item~\ref{tech:angle} of Lemma~\ref{lm:tech} immediately implies item~\ref{struct:angle} of this corollary and the part of item~\ref{struct:tang-comp} of this corollary, namely, the claim that every one-sided tangent of $\Sigma$ at $x$ is tangent to exactly one connected component of $\Sigma \setminus \{x\}$.

Due to item~\ref{tech:nice} of Lemma~\ref{lm:tech}, there is a sequence of $(\Sigma, x)$-nice $\varepsilon_k \rightarrow 0$ such that $\varepsilon_{k + 1} = (1 / 2 + o_{\varepsilon_k}(1)) \cdot \varepsilon_k$ for every $k$.
In view of item~\ref{tech:len} of Lemma~\ref{lm:tech},
\begin{equation}\label{eq:krays}
\H(\Sigma \cap B_{\varepsilon_k}(x) \setminus B_{\varepsilon_{k + 1}}(x) ) = (\varepsilon_k - \varepsilon_{k + 1}) \cdot ( \oordball_x\Sigma + o_{\varepsilon_k}(1) ).
\end{equation}
Let $t \defeq \oordball_x \Sigma$, and let $\{ A^k_1, \dots, A^k_t\} \defeq \Sigma \cap \partial B_{\varepsilon_k}(x) $.
Then, due to~\eqref{eq:krays}, up to a re-enumeration of $A^k_i$ for each fixed $k$, we have $\angle A^k_ixA^{k + 1}_i=o_{\varepsilon_k}(1)$ for $1 \leq i \leq t$.
Then there exist points $A_1, \dots, A_t$ such that $(A^k_ix] \rightarrow (A_ix]$ for $1 \leq i \leq t$.
Due to item~\ref{tech:angle} of Lemma~\ref{lm:tech}, the inequality $\angle A_ixA_j \geq 2 \pi / 3$ holds for any $i \neq j$.

Note that points $A^k_1, \dots, A^k_t$ belong to $k$ different connected components of $\Sigma \setminus \{x\}$ for any sufficiently big $k$.
Indeed, assume that for an arbitrarily big $k$ there exist $i_1 \neq i_2$ such that $A^k_{i_1}$ and $A^k_{i_2}$ belong to the same component.
Then the set 
\[
\Sigma \setminus B_{\varepsilon_k}(x) \cup S^*(x, \varepsilon_k) \cup \bigcup_{i \neq i_2} [A^k_ix],
\]
where $S^*(x, \varepsilon_k)$ is from Corollary~\ref{main_c}, is connected, has the energy not greater than $\Sigma$ and is shorter than $\Sigma$, which is not possible.

Thus to show that $(A_1 x], \dots, (A_t x]$ are $t$ one-sided tangents each tangent to its own connected component of $\Sigma \setminus \{x\}$, it is sufficient to show that for any $z\in \Sigma \setminus \{x\}$ we have
\[
\min_{1 \leq i \leq t}\angle zxA_i=o_{|zx|}(1).
\]
Assume that is not so, i.e. there exists $\gamma > 0$ and a point $z \in \Sigma \setminus \{x\}$ such that $|zx|$ is arbitrarily small and
\[
\min_{1 \leq i \leq t} \angle zxA_i \geq \gamma.
\]

Let $\varepsilon'$ be the infimum of $(\Sigma, x)$-nice $\varepsilon \geq |zx|$.
In view of item~\ref{tech:nice} of Lemma~\ref{lm:tech}, we have $|zx| = \varepsilon' \cdot (1 - o_{\varepsilon'}(1))$.
Due to item~\ref{tech:len} of Lemma~\ref{lm:tech}, $\H(\Sigma \cap B_{\varepsilon'}(x) \setminus B_{|zx|}(x) ) = o_{\varepsilon'}(\varepsilon')$, so there is $z' \in \Sigma \cap \partial B_{\varepsilon'}(x)$ such that $|zz'| = o_{\varepsilon'}(\varepsilon')$, i.e.
\[
\min_{1 \leq i \leq t} \angle z'xA_i \geq \gamma - o_{\varepsilon'}(1).
\]
Let $k$ be such that $\varepsilon_{k + 1} \leq \varepsilon' < \varepsilon_k$. Note that for any $1 \leq i \leq t$ we have
\begin{align*}
|A^k_iz'| + |z'A^{k + 1}_i| &\geq \sqrt{(\sin \gamma - o_{\varepsilon'}(1))^2 \cdot (\varepsilon')^2 + (\varepsilon_k - \varepsilon')^2} + \sqrt{(\sin \gamma - o_{\varepsilon'}(1))^2 \cdot \varepsilon_{k + 1}^2 + (\varepsilon' - \varepsilon_{k + 1})^2}\\
&\geq (\varepsilon_k - \varepsilon_{k + 1}) \cdot \left( \sqrt{\sin^2 \gamma + 1} + o_{\varepsilon_k}(1) \right) \geq (\varepsilon_k - \varepsilon_{k + 1}) \cdot \left( 1 + \frac{\sin^2 \gamma}{3} + o_{\varepsilon_k}(1) \right).
\end{align*}
Then, due to item~\ref{tech:len} of Lemma~\ref{lm:tech},
\begin{align*}
(t + o_{\varepsilon_k}(1)) \cdot \varepsilon_k &= \H( \Sigma \cap B_{\varepsilon_k}(x) )\\
&= \H( \Sigma \cap B_{\varepsilon_{k + 1}}(x) ) + \H( \Sigma \cap B_{\varepsilon'}(x) \setminus B_{\varepsilon_{k + 1}}(x)) + \H( \Sigma \cap B_{\varepsilon_k}(x) \setminus B_{\varepsilon'}(x))\\
&\geq t \cdot \varepsilon_{k + 1} + (t - 1) \cdot (\varepsilon_k - \varepsilon_{k + 1}) + (\varepsilon_k - \varepsilon_{k + 1}) \cdot  \left( 1 + \frac{\sin^2 \gamma}{3} + o_{\varepsilon_k}(1) \right)\\
&\geq \left( t + \frac{\sin^2 \gamma}{6} + o_{\varepsilon_k}(1) \right) \cdot \varepsilon_k,
\end{align*}
which is a contradiction.

We have shown that there are $t = \oordball_x \Sigma \geq \oord_x \Sigma$ one-sided tangents of $\Sigma$ at $x$.
Each of them is tangent to its own connected component of $\Sigma \setminus \{x\}$.
Since $\oord_x \Sigma_1 > 0$ for any connected component $\Sigma_1$ of $\Sigma \setminus \{x\}$, it follows that $\oord_x \Sigma = \oordball_x \Sigma$, and that there are no more connected components, i.e. each connected component of $\Sigma \setminus \{x\}$ indeed has a one-sided tangent.
That concludes the proof of items~\ref{struct:ord} and~\ref{struct:tang-comp} of this corollary.

Finally, item~\ref{struct:comp} follows from previous items of this corollary, Remark~\ref{geq_ord} and item~\ref{tech:len} of Lemma~\ref{lm:tech}.
And item~\ref{struct:loops} follows from item~\ref{struct:comp} of this corollary as for a point $x$ of a loop $C \subset \Sigma$ and for the connected component $\Sigma_1$ of $\Sigma \setminus \{x\}$ such that $\Sigma_1 \supset C\setminus \{x\}$ the order $\oord_x\Sigma_1\geq 2$ which contradicts to item~\ref{struct:comp}.
\end{proof}

\subsection{Technical results}\label{subsec:tech}

In view of Corollary~\ref{cor:struct}, from this point on we are going to make no distinction between $\oord_x \Sigma$ and $\oordball_x \Sigma$, and will always write $\oord_x \Sigma$.

\begin{lemma}\label{napr}
Let $\Sigma$ be a local minimizer for a compact set $M \subset \mathbb{R}^n$ and $r > 0$, let $x \in \G_\Sigma$ be an energetic point and let $\{ (a_ix] \}_{i=1}^{\oord_x \Sigma}$ be the set of all one-sided tangents of $\Sigma$ at $x$. Let
\[
N(x) \defeq \left\{ \vv{z} \in \mathbb{R}^n \ \bigg |\ \vv{z} = \sum_{i=1}^n \alpha_i \cdot \vv{xy_i(x)},\ y_i(x) \text{ corresponds to $x$},\ \alpha_i \geq 0 \right\},
\]
i.e. $N(x)$ is the conical hull of all vectors starting at $x$ and ending at some point, corresponding to $x$.
Then 
\[
\vv{s}(x) \defeq \sum_{i = 1}^{\oord_x \Sigma} \frac{\vv{a_ix}}{|a_ix|} \in N(x).
\]
\end{lemma}
\begin{proof}
Assume the contrary, then $\{\vv{s}(x)\}$ and $N(x)$ are two disjoint nonempty closed convex sets. Moreover, the set $\{\vv{s}(x)\}$ is compact, and $N(x)$ is a cone, so there exists a hyperspace separating $\{\vv{s}(x)\}$ from $N(x)$, i.e. a unit vector $\vv{h} \in \mathbb{R}^n$ such that $\langle \vv{s}(x),\vv{h} \rangle = -\rho < 0$ and $\langle \vv{n} ,\vv{h} \rangle \geq 0$ for any $\vv{n} \in N(x)$.

Let $\varepsilon$ be $(\Sigma, x)$-nice, and let $\{ A_i \}_{i=1}^{\oord_x \Sigma} \defeq \Sigma \cap \partial B_\varepsilon(x)$ be such that $\angle a_ixA_i = o_\varepsilon(1)$.
Let $\vv{S} \defeq \sum_i \vv{A_i x}$, then $\langle \vv{S}, \vv{h} \rangle = -\rho \varepsilon + o_\varepsilon(\varepsilon)$.
By Lemma~\ref{Ooo}, where exists $\delta = o_\varepsilon(\varepsilon)$ such that $|ux| \leq r + \delta$ for any $u \in U(x, \varepsilon)$.
Let $x' \defeq x + k\delta \vv{h}$ for some constant $k > 0$ which we will choose later. Then
\[
\sum_i |A_ix'| = \sum_i |\vv{A_ix} + \vv{xx'}| = \sum_i \sqrt{\varepsilon^2 + k^2\delta^2 + 2 k\delta \langle \vv{A_ix}, \vv{h} \rangle} = \sum_i \varepsilon \cdot \sqrt{ 1 + (k\delta / \varepsilon)^2 + 2 (k\delta / \varepsilon^2) \cdot \langle \vv{A_ix}, \vv{h} \rangle } =
\]
\[
= \sum_i \varepsilon \cdot \left ( 1 + (k\delta / \varepsilon^2) \langle \vv{A_ix}, \vv{h} \rangle + o_\varepsilon(\delta / \varepsilon) \right ) = \oord_x \Sigma \cdot \varepsilon + \langle \vv{S}, \vv{h} \rangle \cdot (k\delta / \varepsilon) + o_\varepsilon(\delta) = \oord_x \Sigma \cdot \varepsilon - \rho k\delta + o_\varepsilon(\delta).
\]
Let $u \in U(x, \varepsilon)$. By Remark~\ref{rem_cC}, there exists a corresponding point $y(x)$ such that $\angle uxy(x) = o_\varepsilon(1)$, so $\langle \vv{xu}, \vv{h} \rangle \geq -o_\varepsilon(1)$. Then 
\[
|x'u| = |\vv{x'x} + \vv{xu}| \leq \sqrt{k^2\delta^2 + ( r + \delta )^2 - 2 k\delta \langle \vv{xu}, \vv{h} \rangle} \leq \sqrt{ ( r + \delta )^2 + o_\varepsilon(\delta) } = r + \delta + o_\varepsilon(\delta).
\]
Then, by Lemma~\ref{cube}, there exists a connected set $S \ni x'$ of length $c \cdot (\delta + o_\varepsilon(\delta))$ for some fixed $c > 0$ such that $U(x, \varepsilon) \subset \overline{B_{r + \delta + o_\varepsilon(\delta)}(x')} \subset \overline{B_r(S)}$. Then the set
\[
\Sigma \setminus B_\varepsilon(x) \cup \bigcup_{i=1}^{\oord_x \Sigma} [A_ix'] \cup S
\]
is connected, is shorter than $\Sigma$ if $k > c/\rho$ and $\varepsilon$ is small enough, and has the energy not larger than $\Sigma$.
It is not possible, since $\Sigma$ is a local minimizer.
\end{proof}

\begin{proposition}\label{prop:alm-orth}
Let $\Sigma$ be a local minimizer for a compact set $M \subset \mathbb{R}^n$ and $r > 0$.
Let $x \in \E_\Sigma$ be a non-isolated energetic point, and let $\Sigma_1$ be a connected component of $\Sigma \setminus \{x\}$ with one-sided tangent $(ax]$.
\begin{enumerate}
\item \label{alm-orth:item1} Let $\Sigma_1$ contain a sequence of energetic points $x_k \in \Sigma_1 \cap \G_\Sigma$ converging to $x$: $x_k \rightarrow x$.
Let a sequence of their corresponding points $y(x_k)$ converge to a point $y$: $y(x_k) \rightarrow y$.
Then $y = y(x)$ corresponds to $x$, and $\vv{xy}$ is orthogonal to $\vv{ax}$.
\item \label{alm-orth:item2} Let $\bar x \in \Sigma_1 \cap \G_\Sigma$ be an energetic point, and let $y(\bar x)$ be any of its corresponding points.
Then there exists a point $y = y(x)$ corresponding to $x$ such that vector $\vv{\bar xy(\bar x)}$ is $o_{|\bar xx|}(1)$-parallel to $\vv{xy}$ and vector $\vv{xy}$ is orthogonal to $\vv{ax}$.
\item \label{alm-orth:ord1} Let $\bar x \in \Sigma_1$ be a point with $\oord_{\bar x} \Sigma = 1$.
Note that, due to Remark~\ref{rem:non-energ}, $\bar x$ is energetic: $\bar x \in \G_\Sigma$.
Then the one-sided tangent $(\bar a\bar x]$ of $\Sigma$ at $\bar x$ is $o_{|\bar xx|}(1)$-orthogonal to $(ax)$.
\item \label{alm-orth:ord2} Let $\bar x \in \Sigma_1$ be a point with $\oord_{\bar x} \Sigma = 2$, and let $(\bar a\bar x]$ and $(\bar b\bar x]$ be one-sided tangents of $\Sigma$ at $\bar x$.
Then $| \angle ((\bar a\bar x], (ax)) - \angle ((\bar b \bar x], (ax)) | = o_{|\bar xx|}(1)$.
\end{enumerate}
\end{proposition}
\begin{proof}
Let us start with a claim slightly weaker than item~\ref{alm-orth:item2}.
\begin{enumerate}[(a)]
\item\label{alm-orth:item0} \textit{Let $\bar x \in \Sigma_1 \cap \G_\Sigma$ be an energetic point, and let $y(\bar x)$ be any of its corresponding points.
Then $\vv{\bar xy(\bar x)}$ is $o_{|\bar xx|}(1)$-orthogonal to $\vv{ax}$.}

Indeed, note that $y(\bar x) \in U(x, |\bar xx|)$, so, by Remark~\ref{rem_cC}, $\angle y(x)xy(\bar x) = o_{|\bar xx|}(1)$ for some corresponding point $y(x)$.
Then, on the one hand, $\Sigma \cap B_r(y(x)) = \emptyset$, so
\[
\angle \bar xxy(x) \geq \pi / 2 - o_{|\bar xx|}(1), \text{ and } \angle \bar xxy(\bar x) \geq \angle \bar xxy(x) - \angle y(\bar x)xy(x) \geq \pi / 2 - o_{|\bar xx|}(1).
\]
On the other hand, $\Sigma \cap B_r(y(\bar x)) = \emptyset$, so $\angle x\bar xy(\bar x) \geq \pi / 2 - o_{|\bar xx|}(1)$.
Then, considering the sum of angles in triangle $\triangle \bar xxy(\bar x)$, we have $\angle x\bar xy(\bar x) = \pi / 2 + o_{|\bar xx|}(1)$, i.e. $\vv{\bar xy(\bar x)}$ is $o_{|\bar xx|}(1)$-orthogonal to $\vv{\bar xx}$.
Finally, $\angle \bar xxa = o_{|\bar xx|}(1)$, so $\vv{\bar xy(\bar x)}$ is $o_{|\bar xx|}(1)$-orthogonal to $\vv{ax}$.
\end{enumerate}
Now let us turn to the claims of the Proposition.
\begin{enumerate}
\item Due to Remark~\ref{emptyball}, $y$ corresponds to $x$.
In view of~\ref{alm-orth:item0}, $\vv{x_ky(x_k)}$ is $o_{|x_kx|}(1)$-orthogonal to $\vv{ax}$, so $\vv{xy}$ is orthogonal to $\vv{ax}$.
\item Assume the contrary, then there exist a constant $c > 0$ and sequences of energetic points $x_k \in \Sigma_1 \cap \G_\Sigma$ and their corresponding points $y(x_k)$ such that $x_k \rightarrow x$ and $\left |\angle (\vv{x_ky(x_k)}, \vv{xy}) \right| \geq c$ for any $k$ and any corresponding to $x$ point $y = y(x)$ with $\vv{xy}$ orthogonal to $\vv{ax}$.
One can choose a converging subsequence $y(x_{k_i}) \rightarrow y$.
Then, due to item~\ref{alm-orth:item1}, $y$ corresponds to $x$, and $\vv{xy}$ is orthogonal to $\vv{ax}$.
But $x_{k_i} \rightarrow x$ and $y(x_{k_i}) \rightarrow y$, so $\angle \left (\vv{x_{k_i}y(x_{k_i})}, \vv{xy} \right) \rightarrow 0$, which is a contradiction.
\end{enumerate}
Item~\ref{alm-orth:ord1} and item~\ref{alm-orth:ord2} for energetic $\bar x$ immediately follow from item~\ref{alm-orth:item2} and Lemma~\ref{napr}.
Item~\ref{alm-orth:ord2} for non-energetic $\bar x$ follows from Remark~\ref{rem:non-energ}.
\end{proof}

\begin{definition}
For any set $\Gamma \subset \mathbb{R}^n$ such that there is a unique path between any pair of points in $\Gamma$ and for any $a, b \in \Gamma$ we will denote by $\Gamma(a, b)$ the unique path in $\Gamma$ connecting $a$ and $b$.
\end{definition}

\begin{definition}
For any points $x_1, \dots, x_k \in \mathbb{R}^n$ such that any two segments from $[x_1x_2], \dots, [x_{k - 1}x_k]$ intersect only, possibly, at endpoints,
we will denote by $x_1{-}\cdots{-}x_k$ a polygonal chain connecting $x_1, \dots, x_k$, i.e. a union of segments $[x_1x_2], \dots, [x_{k - 1}x_k]$.
\end{definition}

\begin{lemma}\label{lm:seqseg}
Let $\Sigma$ be a local minimizer for a compact set $M \subset \mathbb{R}^n$ and $r > 0$.
Let $x \in \Sigma$, and let $\Sigma_1$ be a connected component of $\Sigma \setminus \{x\}$ with one-sided tangent $(ax]$.
Let $x_i \in \Sigma_1$ be a sequence such that $x_i \rightarrow x$ and the path $\Sigma(x_1, x_i)$ coincides with the polygonal chain $x_1{-}\cdots{-}x_i$
for any $i > 1$.
Then there exists an infinite subsequence $x_{i_k}$ such that $\angle((x_{i_k}x_{i_k + 1}), (ax)) \rightarrow 0$.
\end{lemma}
\begin{proof}
Assume the contrary, i.e. there exists a positive $c < \pi / 2$ such that $|\angle((x_ix_{i + 1}), (ax))| \geq c$ for every sufficiently big $i$.
Denote by $\bar x_i$ the orthogonal projection of $x_i$ onto $(ax)$, then $|\bar x_i\bar x_{i + 1}| \leq |x_ix_{i + 1}| \cdot \cos c$ for every sufficiently big $i$.
By the definition of a one-sided tangent, $\angle x_ixa \rightarrow 0$, so $|\bar x_ix| = |x_ix| \cdot (1 + o_{|x_ix|}(1))$.
Due to item~\ref{tech:nice} of Lemma~\ref{lm:tech}, for every sufficiently big $i$ there exists a $(\Sigma, x)$-nice $\delta_i \geq |x_ix|$ such that $\delta_i = |x_ix| \cdot (1 + o_{|x_ix|}(1))$.
Note that $\Sigma(x_i, x) \subset B_{\delta_i}(x)$ if $i$ is big enough.
Then for every sufficiently  big $i$ we have
\[
\H(\Sigma_1 \cap B_{\delta_i}(x)) \geq \H(\Sigma(x_i, x)) = \sum_{j = i}^{\infty} |x_jx_{j + 1}| \geq \sum_{j = i}^{\infty} \frac{|\bar x_j\bar x_{j + 1}|}{\cos c} \geq \frac{|\bar x_ix|}{\cos c} = \delta_i \cdot \left( \frac{1}{\cos c} + o_{|x_ix|}(1) \right),
\]
which contradicts item~\ref{struct:comp} of Corollary~\ref{cor:struct}.
\end{proof}

The following proposition is an analogue of Remark~\ref{rem:non-energ} for isolated energetic points.
Even though it is true in general case ($\mathbb{R}^n$ for $n > 2$), in view of Example~\ref{rem:all-energ}, we are able to use it in a meaningful way only in planar case.

\begin{proposition}\label{discr_str}
Let $\Sigma$ be a local minimizer for a compact set $M \subset \mathbb{R}^n$ and $r > 0$.
For any isolated energetic point $x\in \X_\Sigma$ there exists $\rho = \rho(x) > 0$ such that the set $\Sigma \cap \overline{B_\rho(x)}$ is a union of $\oord_x \Sigma$ segments with one end at $x$ and another end at $\partial B_\rho(x)$, with pairwise angles between them at least $2 \pi / 3$.
In view of Proposition~\ref{3vectors}, these segments are coplanar.
\end{proposition}
\begin{proof}
Let $\varepsilon$ be $(\Sigma, x)$-nice and such that $\G_\Sigma \cap \overline{B_\varepsilon(x)} = \{x\}$.
There exists a sequence of $(\Sigma, x)$-nice $\varepsilon_i < \varepsilon$ such that $\varepsilon_i \rightarrow 0$.
For each $i$ define the set
\[
\Sigma_i \defeq \Sigma \cap \overline{B_\varepsilon(x)} \setminus B_{\varepsilon_i}(x).
\]
We can assume that $\varepsilon$ is small enough that $\Sigma$ is the shortest set among all sets $\Sigma'$ with $\diam(\Sigma \triangle \Sigma') \leq \varepsilon$ and $\F_M(\Sigma') \leq r$.
Thus a closed set $\Sigma_i \subset \S_\Sigma$ is a Steiner forest.
Note that $\Sigma_i$ has $\oord_x \Sigma$ connected components, one corresponding to each connected component of $\Sigma \setminus \{x\}$.
Moreover, each component of $\Sigma_i$ is a Steiner tree for the set of two vertices (one from $\partial B_\varepsilon(x)$ and another from $\partial B_{\varepsilon_i}(x)$), i.e. a segment.
It follows that the set $\Sigma \cap \overline{B_\varepsilon(x)} = \{x\} \cup \bigcup_{i = 1}^{\infty} \Sigma_i$ is also a union of $\oord_x \Sigma$ segments with one end at $x$ and another at $\partial B_\varepsilon(x)$.
In view of item~\ref{struct:angle} of Corollary~\ref{cor:struct}, the pairwise angles between these segments are at least $2 \pi / 3$.
\end{proof}


\section{Planar results}
\label{planarcase}

Note that an energetic point $x \in \G_\Sigma$ in $\mathbb{R}^2$ cannot have $\oord_x \Sigma > 2$ (due to item~\ref{struct:angle} of Corollary~\ref{cor:struct}), but it is possible in $\mathbb{R}^n$ when $n > 2$ (see Example~\ref{rem:all-energ}).

\begin{lemma}\label{lm:ord2parallel}
Let $\Sigma$ be a local minimizer for a compact set $M \subset \mathbb{R}^2$ and $r > 0$, let $x \in \E_\Sigma$ be a non-isolated energetic point, and let $\Sigma_1$ be a connected component of $\Sigma \setminus \{x\}$ with one-sided tangent $(ax]$.
Then there exists $\rho = \rho(x) > 0$ such that, if $\bar x \in \Sigma_1 \cap B_\rho(x) \cap \E_\Sigma$ is a non-isolated energetic point, then $\oord_{\bar x} \Sigma = 2$, and every one-sided tangent $(\bar a\bar x]$ of $\Sigma$ at $\bar x$ is $o_{|\bar xx|}(1)$-parallel to $(ax)$.
\end{lemma}
\begin{proof}
Let $x_k \rightarrow \bar x$ and $y(x_k)$ be sequences satisfying conditions of item~\ref{alm-orth:item1} of Proposition~\ref{prop:alm-orth}; let $(\bar a\bar x]$ be the one-sided tangent to the connected component of $\Sigma_1 \setminus \{\bar x\}$ containing $x_k$, and let $y(x_k) \rightarrow y(\bar x)$.
Then $\vv{\bar xy(\bar x)}$ is orthogonal to $\vv{\bar a\bar x}$.
On the other hand, in view of item~\ref{alm-orth:item2} of Proposition~\ref{prop:alm-orth}, $\vv{\bar xy(\bar x)}$ is $o_{|\bar xx|}(1)$-orthogonal to $\vv{ax}$.
It follows that $(\bar a\bar x]$ is $o_{|\bar xx|}(1)$-parallel to $(ax)$.

If $\oord_{\bar x} \Sigma = 1$, then, in view of item~\ref{alm-orth:ord1} of Proposition~\ref{prop:alm-orth}, $(\bar a\bar x]$ would be simultaneously $o_{|\bar xx|}(1)$-parallel and $o_{|\bar xx|}(1)$-orthogonal to $(ax)$, which is impossible when $\rho \geq |\bar xx|$ is sufficiently small.
On the other hand, $\oord_{\bar x} \Sigma \neq 3$ since a branching point is always non-energetic in $\mathbb{R}^2$, so $\oord_{\bar x} \Sigma = 2$.

Finally, the second one-sided tangent $(\bar b\bar x]$ of $\Sigma$ at $\bar x$ is also $o_{|\bar xx|}(1)$-parallel to $(ax)$ due to item~\ref{alm-orth:ord2} of Proposition~\ref{prop:alm-orth}.
\end{proof}

\begin{proposition}\label{brrrr}
Let $\Sigma$ be a local minimizer for a compact set $M \subset \mathbb{R}^2$ and $r > 0$, let $x \in \E_\Sigma$ be a non-isolated energetic point, and let $\Sigma_1$ be a connected component of $\Sigma \setminus \{x\}$ with one-sided tangent $(ax]$.
Then for any polygonal chain $x_1{-}x_2{-}x_3{-}x_4{-}x_5 \subset \Sigma_1$, where $x_i \in \X_\Sigma$, $\oord_{x_i} \Sigma = 2$ for $1 \leq i \leq 5$, one has $\max_{1 \leq i \leq 4} \angle((x_ix_{i + 1}) , (ax)) = o_\rho(1)$, where $\rho \defeq \max_{1 \leq i \leq 5} |x_ix|$.
\end{proposition}
\begin{proof}
Assume the contrary, i.e. there exists $c > 0$ such that for any $\rho > 0$ there exists a polygonal chain described above satisfying $|x_ix| \leq \rho$ for $1 \leq i \leq 5$ and $\angle((x_jx_{j + 1}), (ax)) \geq c$ for some $1 \leq j \leq 4$.
Let $y_i = y_i(x_i)$ be an arbitrary corresponding to $x_i$ point for $1 \leq i \leq 5$. Due to item~\ref{alm-orth:item2} of Proposition~\ref{prop:alm-orth},
\[
\angle(\vv{x_iy_i}, \vv{ax}) = \pi / 2 + o_\rho(1) \text{ for } 1 \leq i \leq 5.
\]
In view of item~\ref{alm-orth:ord2} of Proposition~\ref{prop:alm-orth}, $|\angle((x_ix_{i + 1}), (ax)) - \angle((x_ix_{i - 1}), (ax))| = o_\rho(1)$ for $2 \leq i \leq 4$, so
\begin{equation}\label{eq:seg-ax-angle}
\angle((x_ix_{i + 1}), (ax)) = \alpha + o_\rho(1) \text{ for } 1 \leq i \leq 4 \text{ and some } \alpha \geq c.
\end{equation}

By the definition of a corresponding point, $\angle(\vv{x_iy_i}, \vv{x_ix_{i \pm 1}}) \geq \pi / 2$ for $1 \leq i \leq 5$ whenever $x_{i \pm 1}$ is defined, so
\[
\angle(\vv{x_iy_i}, \vv{x_ix_{i \pm 1}}) = \pi / 2 + \alpha + o_\rho(1) \text{ for } 1 \leq i \leq 5 \text{ whenever } x_{i \pm 1} \text{ is defined.}
\]
It follows that
\begin{equation}\label{eq:segs-diff}
\angle(\vv{x_ix_{i - 1}}, \vv{x_ix_{i + 1}}) = \pi - 2\alpha + o_\rho(1) \text{ for } 2 \leq i \leq 4
\end{equation}
and
\[
\angle(\vv{x_iy_i}, \vv{x_{i + 1}y_{i + 1}}) = \pi + o_\rho(1) \text{ for } 1 \leq i \leq 4.
\]

For $2 \leq i \leq 4$, by orthogonally projecting each summand in the equation 
\[
\vv{y_{i - 1}x_{i + 1}} = \vv{y_{i - 1}x_{i - 1}} + \vv{x_{i - 1}x_i} + \vv{x_ix_{i + 1}}
\]
onto $(y_{i - 1}x_{i + 1})$, we get 
\[
r \leq |y_{i - 1}x_{i + 1}| \leq r + |x_{i - 1}x_i| (1 + o_\rho(1)) \sin\alpha - |x_ix_{i + 1}| (1 + o_\rho(1)) \sin\alpha,
\]
i.e. $|x_{i - 1}x_i| / |x_ix_{i + 1}| \geq 1 + o_\rho(1)$.
Applying the same arguments to $\vv{y_{i + 1}x_{i - 1}}$, we get $|x_ix_{i + 1}| / |x_{i - 1}x_i| \geq 1 + o_\rho(1)$, so
\begin{equation}\label{eq:same-length}
|x_ix_{i + 1}| = l \cdot (1 + o_\rho(1)) \text{ for } 1 \leq i < 5 \text{ and some } 0 < l \leq 2\rho.
\end{equation}

It follows from~\eqref{eq:seg-ax-angle},~\eqref{eq:segs-diff} and~\ref{eq:same-length} that $\angle ((x_2x_4), (x_2x_3)) = \alpha + o_\rho(1)$ and $\angle ((x_2x_4), (x_3x_4)) = \alpha + o_\rho(1)$, so
\[
|x_2x_3| + |x_3x_4| - |x_2x_4| = 2 (1 - \cos\alpha)(1 + o_\rho(1)) \cdot l.
\]

Identities~\eqref{eq:seg-ax-angle},~\eqref{eq:segs-diff} and~\ref{eq:same-length} also imply that $\angle((x_1x_5), (x_1x_3)) = o_\rho(1)$, so $\delta \defeq \dist([x_1x_5], x_3) = o_\rho(l)$.
It follows that $( [x_2x_3] \cup [x_3x_4] ) \subset \conv(\{x_1, x_2, x_4, x_5\}) \cup B_\delta([x_1x_5])$, i.e.
\[
\overline{B_r([x_2x_3] \cup [x_3x_4])} \subset \overline{B_r([x_1x_2] \cup [x_2x_4] \cup [x_4x_5])} \cup \overline{B_{r + \delta}([x_1x_5])}.
\]

Note that $\overline{B_{r + \delta}([x_1x_5])} \subset \overline{B_{r + \delta'}(\{x_1, x_5\})}$ for some $\delta' = o_\rho(l)$, since
\[
\sqrt{(r + \delta)^2 + \left(\frac{|x_1x_5|}{2}\right)^2} = \sqrt{r^2 + o_\rho(l)} = r + o_\rho(l).
\]
Let the sets $S_j \defeq S(x_j, \delta')$, $\H(S_j) = c\delta' = o_\rho(l)$, $j \in \{1, 5\}$, be from Lemma~\ref{cube}.
We see that the set
\[
\Sigma \setminus [x_2x_3] \setminus [x_3x_4] \cup [x_2x_4] \cup S_1 \cup S_5
\]
is connected and has energy not larger than $\Sigma$.
It is also shorter than $\Sigma$ if $\rho$ is small enough, which contradicts the definition of a local minimizer.
\end{proof}

\begin{theorem}\label{thm:plane}
Let $\Sigma$ be a local minimizer for a compact set $M \subset \mathbb{R}^2$ and $r > 0$.
There is a finite number of branching points in $\Sigma$.
\end{theorem}
\begin{proof}
Assume the contrary, then the set of branching points is infinite and thus has a limit point.
Let $x \in \Sigma$ be such a limit point, and let $\Sigma'$ be a connected component of $\Sigma \setminus \{x\}$ which contains an infinite sequence of branching points converging to $x$.
Let $(ax]$ be the unique one-sided tangent of $\Sigma'$ at $x$.
Note that, due to Remark~\ref{rem:non-energ} and Proposition~\ref{discr_str}, $x$ is a non-isolated energetic point: $x \in \E_\Sigma$.

Given a $(\Sigma, x)$-nice $\varepsilon > 0$, let $x'_\varepsilon$ be the unique point in $\Sigma' \cap \partial B_\varepsilon(x)$.
Also, for any branching point $z \in \Sigma' \cap B_\varepsilon(x)$ let us denote by $\Sigma_z$ the unique connected component of $\Sigma \setminus \{z\}$ not containing either of $x, x'_\varepsilon$.
Let us study the structure of $\Sigma_z$ when $\varepsilon$ is small enough.

\begin{enumerate}[(a)]
\item \label{plane:uniquey} \textit{There exists $\rho = \rho(x) > 0$ such that for any $(\Sigma, x)$-nice $\varepsilon < \rho$ and any branching point $z \in \Sigma' \cap B_\varepsilon(x)$ there exists a corresponding to $x$ point $y = y(x)$ with $\vv{xy(x)}$ orthogonal to $\vv{ax}$ such that for every energetic point $\bar x \in \Sigma_z \cap \G_\Sigma$ and its corresponding point $y(\bar x)$ vector $\vv{\bar xy(\bar x)}$ is $o_{\varepsilon}(1)$-parallel to $\vv{xy(x)}$.
Moreover,
\begin{itemize}
\item if $\oord_{\bar x} \Sigma = 1$, and $(\bar a\bar x]$ is the one-sided tangent of $\Sigma$ at $\bar x$, then $\vv{\bar a\bar x}$ is $o_{\varepsilon}(1)$-parallel to $\vv{xy(x)}$;
\item if $\oord_{\bar x} \Sigma = 2$, $(\bar a\bar x]$ and $(\bar b\bar x]$ are one-sided tangents of $\Sigma$ at $\bar x$, then $| \angle (\vv{\bar a\bar x}, \vv{xy(x)}) - \angle (\vv{\bar b \bar x}, \vv{xy(x)}) | = o_{\varepsilon}(1)$.
\end{itemize}
}

Assume the contrary, then, in view of item~\ref{alm-orth:item2} of Proposition~\ref{prop:alm-orth}, there exist two corresponding to $x$ points $y_1(x)$ and $y_2(x)$ with vectors $\vv{xy_1(x)}$ and $\vv{xy_2(x)}$ orthogonal to $\vv{ax}$ (so $\vv{xy_1(x)} = -\vv{xy_2(x)}$) such that for $i \in \{1, 2\}$ there exists an energetic point $\bar x_i \in \Sigma_z \cap \G_\Sigma$ with a corresponding point $y(\bar x_i)$ such that vector $\vv{\bar x_iy(\bar x_i)}$ is $o_{\varepsilon}(1)$-parallel to $\vv{xy_j(x)}$.

By the definition of a one-sided tangent, any point $p \in \Sigma' \cap B_\varepsilon(x)$ satisfies $\angle pxa < \delta$ for some $\delta = o_\varepsilon(1)$.
Let $Q \defeq \{ p \in B_\varepsilon(x) \ |\ \angle pxa \leq \delta\}$, then $\Sigma' \cap B_\varepsilon(x) \subset Q$.
Note that $\Sigma(x'_\varepsilon, x)$ splits the set $Q$ into two simply connected sets; denote by $A$ one of these sets which contains both $\bar x_1$ and $\bar x_2$.
For $i \in \{1, 2\}$ denote $B_i \defeq \{ p \in B_{r + \varepsilon}(x) \setminus B_{r - \varepsilon}(x) \ |\ \angle pxy_i(x) \leq \gamma \}$, where $\gamma \defeq \max_i \angle y(\bar x_i)xy_i(x) = o_{\varepsilon}(1)$.
By definition, $y(\bar x_i) \in B_i$ for $i \in \{1, 2\}$.
Note that it is true for one of $i \in \{1, 2\}$ that any segment connecting a point from $B_i$ with a point from $A$ intersects $\Sigma(x'_\varepsilon, x)$.
It follows that $\Sigma(x'_\varepsilon, x)$ intersects one of segments $[\bar x_1y(\bar x_1)]$, $[\bar x_2y(\bar x_2)]$ at an internal point of a segment, which contradicts the definition of a corresponding point.

\item \label{plane:branch-seg} \textit{There exists $\rho = \rho(x) > 0$ such that for any $(\Sigma, x)$-nice $\varepsilon < \rho$ and for every branching point $z \in \Sigma' \cap B_\varepsilon(x)$ the set $\Sigma_z$ is a segment which is $o_\varepsilon(1)$-perpendicular to $(ax)$.}

Let us consider an arbitrary point $t_1$ of $\Sigma_z$ with $\oord_{t_1} \Sigma = 1$.
We are going to build a sequence of points $t_i \in (\Sigma_z \cup \{z\})$ such that
\begin{enumerate}[i.]
    \item \label{branch-seg:poly} $\Sigma(t_1, t_i) = t_1{-}\cdots{-}t_i$ for every $i > 1$;
    \item \label{branch-seg:angle} $\angle (\vv{t_it_{i - 1}}, \vv{xy(x)}) \leq \pi / 3 + \gamma$ for every $i > 1$, where $y(x)$ is from~\ref{plane:uniquey} and $\gamma = \angle(\vv{t_2t_1}, \vv{xy(x)}) = o_\varepsilon(1)$;
    \item \label{branch-seg:branch} $t_i$ is a branching point for every $i > 1$.
\end{enumerate}
In view of Remark~\ref{rem:non-energ}, $t_1$ is an energetic point: $t_1 \in \G_\Sigma$, then, in view of Lemma~\ref{lm:ord2parallel}, $t_1$ is an isolated energetic point: $t_1 \in \X_\Sigma$.
Let $(a_1t_1]$ be the one-sided tangent of $\Sigma$ at $t_1$, then, due to~\ref{plane:uniquey}, $\gamma \defeq \angle(\vv{a_1t_1}, \vv{xy(x)}) = o_\varepsilon(1)$.
By Proposition~\ref{discr_str}, we have $((a_1t_1] \cap \overline{B_\delta(t_1)}) \subset \Sigma$ for any sufficiently small $\delta > 0$;
choose maximal possible such $\delta$ and define $t_2$ as the unique point in $(a_1t_1] \cap \partial B_\delta(t_1)$.
Items~\ref{branch-seg:poly} and~\ref{branch-seg:angle} for $i = 2$ hold by the construction.
Note that~\ref{branch-seg:angle} and~\ref{plane:uniquey} imply that $t_2$ is non-energetic if $\varepsilon$ is small enough, so, due to Remark~\ref{rem:non-energ}, $\oord_{t_2} \Sigma \in \{2, 3\}$.
Since $\delta$ in the definition of $t_2$ was maximal possible, from the same Remark~\ref{rem:non-energ} it follows that $\oord_{t_2} \Sigma \neq 2$, so~\ref{branch-seg:branch} for $i = 2$ also holds.

Now, suppose we have already chosen points $t_1, \dots, t_j$, $j > 1$, such that conditions~\ref{branch-seg:poly},~\ref{branch-seg:angle},~\ref{branch-seg:branch} hold for $i \leq j$.
If $t_j = z$, we stop.
Otherwise, since $t_j$ is a branching point, there exists a one-sided tangent $(a_jt_j]$ of $\Sigma$ at $t_j$ such that $\angle (\vv{a_jt_j}, \vv{xy(x)}) \leq \pi / 3 + \gamma$.
Let $\delta > 0$ be maximal possible such that $((a_jt_j] \cap \overline{B_\delta(t_j)}) \subset \Sigma$, and choose $t_{j + 1}$ as the unique point in $(a_jt_j] \cap \partial B_\delta(t_j)$.
Again,~\ref{branch-seg:poly} and~\ref{branch-seg:angle} for $i = j + 1$ hold by the construction, and the same arguments as in the case $i = 2$ show~\ref{branch-seg:branch} for $i = j + 1$.

Assume that $t_i \neq z$ for every $i$. Then the infinite sequence of branching points $t_i$ converges to some point, which is non-isolated energetic due to Remark~\ref{rem:non-energ} and Proposition~\ref{discr_str}: $t_i \rightarrow \bar x \in \E_\Sigma \cap (\Sigma_z \cup \{z\})$.
Let $(\bar a\bar x]$ be the one-sided tangent to a connected component of $\Sigma \setminus \{\bar x\}$ which contains the sequence $t_i$.
Due to Lemma~\ref{lm:ord2parallel}, $(\bar a\bar x]$ is $o_{\varepsilon}(1)$-parallel to $(ax)$.
Thus, for every $i > 1$ we have
\[
\angle( (t_it_{i - 1}), (\bar a\bar x)) \geq \angle( (t_it_{i - 1}), (ax)) - o_{\varepsilon}(1) \geq \pi / 2 - \pi / 3 - \gamma - o_{\varepsilon}(1) \geq \pi / 7,
\]
if $\varepsilon$ is small enough.
This contradicts Lemma~\ref{lm:seqseg}.

It follows that $t_i = z$ for some $i > 1$.
Note that in the arguments above there were two possible choices for $t_3$, so, if $t_2 \neq z$, we would have two different paths from $t_1$ to $z$ in $\Sigma$, which is not possible.
Therefore, $t_2 = z$, $\Sigma_z = [t_1z]$, and $\angle (\vv{zt_1}, \vv{ax}) = \pi / 2 - \angle(\vv{zt_1}, \vv{xy(x)}) = \pi / 2 - \gamma = \pi / 2 - o_\varepsilon(1)$.
\end{enumerate}
Now let us return to the proof of the theorem.
Let $\varepsilon$ be $(\Sigma, x)$-nice and satisfy $\varepsilon < \rho$, where $\rho$ is from~\ref{plane:branch-seg}.
Recall that $\Sigma'$ contains a sequence of branching points converging to $x$, so there exists a branching point $x_1 \in \Sigma' \cap B_\varepsilon(x)$.
We are going to build a sequence of points $x_i$ such that 
\begin{enumerate}[i*.]
\item \label{seq:poly} $\Sigma(x_1, x_i) = x_1{-}\cdots{-}x_i \subset \Sigma(x_1, x)$ for every $i > 1$;
\item \label{seq:angle} $\angle((x_ix_{i + 1}), (ax)) \geq \pi / 6 - \gamma$ for every $i$, where $\gamma = o_\varepsilon(1)$ does not depend on $i$;
\item \label{seq:point} $x_i$ is either a branching point or an isolated energetic point of order 2 for every $i$.
\end{enumerate}
Suppose $x_1, \dots, x_i$ are already defined.
Let $(a_ix_i]$ be the one-sided tangent of $\Sigma(x_i, x)$ at $x_i$.
Let $\delta > 0$ be maximal possible such that $( (a_ix_i] \cap \overline{B_\delta(x_i)} ) \subset \Sigma(x_i, x)$, and set $x_{i + 1}$ as the unique point in $(a_ix_i] \cap \partial B_\delta(x_i)$.
Item~\ref{seq:poly} follows by the construction.
Note that if $x_i$ is a branching point, then $\angle((x_ix_{i + 1}), (ax)) \geq \pi / 6 - \gamma'$, where $\gamma' = o_\varepsilon(1)$ is from~\ref{plane:branch-seg} and does not depend on $i$.
Otherwise, in view of Proposition~\ref{brrrr}, there is a branching point $x_j$ for some $j \geq i - 4$ if $\varepsilon$ is small enough.
Then we have
\[
\angle((x_ix_{i + 1}), (ax)) \geq \angle((x_jx_{j + 1}), (ax)) - (i - j) \gamma'' \geq \pi / 6 - \gamma' - 4\gamma'' = \pi / 6 - \gamma,
\]
where $\gamma'' = o_\varepsilon(1)$ is from item~\ref{alm-orth:ord2} of Proposition~\ref{prop:alm-orth}, and $\gamma \defeq \gamma' + 4\gamma'' = o_\varepsilon(1)$ does not depend on $i$.
Thus item~\ref{seq:angle} holds.
Finally, to see that item~\ref{seq:point} holds, note that if $x_{i + 1}$ is not a branching point, then it is an energetic point of order 2, which cannot be non-isolated due to Lemma~\ref{lm:ord2parallel}.

Note that conditions~\ref{seq:poly} and~\ref{seq:angle} imply that $x_i \neq x$ for every $i$.
Then the sequence $x_i$ converges to some point $\bar x$: $x_i \rightarrow \bar x$.
Let $(\bar a\bar x]$ be the one-sided tangent to the connected component of $\Sigma \setminus \{\bar x\}$ containing $x_i$.
Note that either $\bar x = x$ and $(\bar a \bar x] = (ax]$, or $(\bar a \bar x]$ is $o_\varepsilon(1)$-parallel to $(ax)$ due to Lemma~\ref{lm:ord2parallel}.
Either way, we have a contradiction between item~\ref{seq:angle} and Lemma~\ref{lm:seqseg}.
\end{proof}

\begin{lemma}\label{zapred2}
Let $\Sigma$ be a local minimizer for a compact set $M \subset \mathbb{R}^2$ and $r > 0$ and let $x \in \Sigma$.
Let $\Sigma_1$ be a connected component of $\Sigma \setminus \{x\}$ with one-sided tangent $(ax]$ and let $\bar x \in \Sigma_1$.
\begin{enumerate}
\item \label{zapred2:weak} For any one-sided tangent $(\bar a\bar x]$ of $\Sigma$ at $\bar x$ the equality $\angle((\bar a\bar x), (ax)) = o_{|\bar xx|}(1)$ holds.
\item Let $(\bar a\bar x]$ be a one-sided tangent at $\bar x$ of any connected component of $\Sigma \setminus \{\bar x\}$ not containing $x$.
Then $\angle((\bar a\bar x], (ax]) = o_{|\bar xx|}(1)$.
\end{enumerate}
\end{lemma}
\begin{proof}
First, note that if $x \notin \E_\Sigma$, then, due to Remark~\ref{rem:non-energ} and Proposition~\ref{discr_str}, there is $\rho > 0$ such that $\Sigma_1 \cap \overline{B_\rho(x)} \subset (ax]$, and claims of the lemma are true.
Thus we can assume that $x \in \E_\Sigma$.
\begin{enumerate}
\item Assume the contrary, i.e. there exists some positive $c < \pi / 2$ and a sequence $x_k \in \Sigma_1$ such that $x_k \rightarrow x$ and $\angle((a_kx_k), (ax)) \geq c$ for any $k$, where $(a_kx_k]$ is some one-sided tangent of $\Sigma$ at $x_k$.
In view of Lemma~\ref{lm:ord2parallel} and Theorem~\ref{thm:plane}, we can assume that there are no branching points and no non-isolated energetic points in the sequence $x_k$.
Then, in view of item~\ref{alm-orth:ord2} of Proposition~\ref{prop:alm-orth}, we can assume that each $(a_kx_k]$ is a one-sided tangent of a connected component of $\Sigma \setminus \{x_k\}$ not containing $x$.

Given $x_k$ for a sufficiently big $k$, we build a polygonal chain $l_1{-}l_2{-}l_3{-}l_4{-}l_5$ as follows: set $l_1 \defeq x_k$ if $x_k \in \X_\Sigma$, otherwise note that $(b_kx_k]$, where $\vv{a_kx_k} = -\vv{b_kx_k}$, is the one-sided tangent of $\Sigma(x_k, x)$ at $x_k$, let $\delta > 0$ be maximal possible such that $( (b_kx_k] \cap \overline{B_\delta(x_k)} ) \subset \Sigma(x_k, x)$, and choose $l_1$ as the unique point in $(b_kx_k] \cap \partial B_\delta(x_k)$.
Note that in the latter case $l_1 \neq x$ as $\angle((x_kl_1), (ax)) = \angle((a_kx_k), (ax)) \geq c$.
In either case, $l_1 \in \X_\Sigma$.

Now, given $l_i \in \X_\Sigma$ for $1 \leq i \leq 4$, let $(bl_i]$ be the one-sided tangent of $\Sigma(l_i, x)$ at $l_i$, let $\delta > 0$ be maximal possible such that $( (bl_i] \cap \overline{B_\delta(l_i)} ) \subset \Sigma(l_i, x)$, and choose $l_{i + 1}$ as the unique point in $(bl_i] \cap \partial B_\delta(l_i)$.
Due to item~\ref{alm-orth:ord2} of Proposition~\ref{prop:alm-orth}, $\angle((l_il_{i + 1}), (ax)) \geq c - o_{|x_kx|}(1)$, so $l_{i + 1} \neq x$ if $k$ is big enough.

We have shown that for any sufficiently big $k$ there exists a polygonal chain $l_1{-}l_2{-}l_3{-}l_4{-}l_5 \subset \Sigma_1$ such that $l_i \in \X_\Sigma$, $\oord_{l_i} \Sigma = 2$ for $1 \leq i \leq 5$, and $\angle((l_il_{i + 1}), (ax)) \geq c - o_{|x_kx|}(1)$ for $1 \leq i \leq 4$.
This is a contradiction to Proposition~\ref{brrrr}.

\item Assume the contrary, then there exists $c > 0$ and $\bar x \in \Sigma_1$ such that $|\bar xx|$ is arbitrarily small and $\angle((\bar a\bar x], (ax]) \geq c$, where $(\bar a\bar x]$ is a one-sided tangent of a connected component of $\Sigma \setminus \{\bar x\}$ not containing $x$.
Due to item~\ref{zapred2:weak}, we have $\angle((\bar a\bar x], (ax]) \geq \pi - o_{|\bar xx|}(1)$.
Let $(\bar b\bar x]$ be the one-sided tangent of $\Sigma(\bar x, x)$ at $\bar x$.
We have $\angle((\bar a\bar x], (\bar b\bar x]) \geq 2 \pi / 3$, so $\angle( (\bar b\bar x], (ax] ) \leq \pi / 2$ if $|\bar xx|$ is small enough.

It follows that $\bar x$ is not the farthest from $x$ point in $\Sigma(\bar x, x)$; let $x' \neq \bar x$ be the farthest from $x$ point in $\Sigma(\bar x, x)$.
Let $(a'_1x']$ and $(a'_2x']$ be one-sided tangents of $\Sigma(\bar x, x')$ and $\Sigma(x', x)$ at $x'$, respectively.
Since $x'$ is farthest from $x$ in $\Sigma(\bar x, x)$, we have $\angle((a'_ix'], (ax]) \geq \pi / 2$ for $i \in \{1, 2\}$.
On the other hand, due to~\ref{zapred2:weak}, $\angle((a'_ix'), (ax)) = o_{|x'x|}(1)$ for $i \in \{1, 2\}$, so $\angle((a'_ix'], (ax]) \geq \pi - o_{|x'x|}(1)$, which contradicts item~\ref{struct:angle} of Corollary~\ref{cor:struct}.
\end{enumerate}
\end{proof}

\section*{Acknowledgments}

Research of Alexey Gordeev is supported by Ministry of Science and Higher Education of the Russian Federation, agreement №075-15-2022-289.
Research of Yana Teplitskaya is supported by ``Native towns'', a social investment program of PJSC ``Gazprom Neft''.
We would like to thank Danila Cherkashin for fruitful discussions and Azat Miftakhov for an inspiring example of doing mathematics in any conditions.
While we were finishing this article Russia invaded Ukraine.
Stand with Ukraine.

\bibliographystyle{amsplain}
\bibliography{main}

\end{document}